\documentclass{svjour3}                     
\smartqed  
\usepackage{color}
\usepackage{amsmath,amsfonts,amssymb,bm,mathtools}
\usepackage{graphicx,psfrag,subfigure}
\usepackage{enumitem}
\usepackage{epstopdf}
\usepackage{latexsym}
\usepackage{mathrsfs}
\usepackage{empheq}
\usepackage{algorithm}
\usepackage{caption}
\usepackage{ulem}

\usepackage{fullpage}

\allowdisplaybreaks

\newtheorem{thm}{Theorem}

\numberwithin{equation}{section}

\begin{document}

\title{A weighted scalar auxiliary variable method for solving gradient flows:
       bridging the nonlinear energy-based and Lagrange multiplier approaches}

\titlerunning{A weighted scalar auxiliary variable method for solving gradient flows}

\author{Qiong-Ao Huang \and Wei Jiang \and Jerry Zhijian Yang \and Cheng Yuan}


\institute{Qiong-Ao Huang\at
School of Mathematics and Statistics, Henan University, Kaifeng 475004, China\\
Center for Applied Mathematics of Henan Province, Henan University, Zhengzhou 450046, China\\
\email{huangqiongao@henu.edu.cn}
\and
Wei Jiang\and Jerry Zhijian Yang\at
School of Mathematics and Statistics, Wuhan University, Wuhan 430072, China\\
Hubei Key Laboratory of Computational Science, Wuhan University, Wuhan 430072, China\\
\email{jiangwei1007@whu.edu.cn; zjyang.math@whu.edu.cn}
\and
Cheng Yuan\at
School of Mathematics and Statistics, and Hubei Key Lab--Math. Sci., Central China Normal University, Wuhan 430079, China\\
Key Laboratory of Nonlinear Analysis \& Applications (Ministry of Education), Central China Normal University, Wuhan 430079, China\\
\email{yuancheng@ccnu.edu.cn}
}

\date{Received: date / Accepted: date}

\maketitle

\begin{abstract}
Two primary scalar auxiliary variable (SAV) approaches are widely applied for simulating gradient flow systems, i.e., the nonlinear energy-based approach and the Lagrange multiplier approach. The former guarantees unconditional energy stability through a modified energy formulation, whereas the latter preserves original energy stability but requires small time steps for numerical solutions. In this paper, we introduce a novel weighted SAV method which integrates these two approaches for the first time. Our method leverages the advantages of both approaches: (i) it ensures the existence of numerical solutions for any time step size with a sufficiently large weight coefficient; (ii) by using a weight coefficient smaller than one, it achieves a discrete energy closer to the original, potentially ensuring stability under mild conditions; and (iii) it maintains consistency in computational cost by utilizing the same time/spatial discretization formulas. We present several theorems and numerical experiments to validate the accuracy, energy stability and superiority of our proposed method.


\end{abstract}

\section{Introduction}
\label{secint}

Gradient flow equations, as an important class of parabolic partial differential equations,
have been widely utiliezd in various fields of science and engineering,
such as wrinkling phenomenon, crystal growth, liquid crystal, image processing, solid-state dewetting,
and geometric flows~\cite{ShenJ19,BaoW23,JiangW23,Huang19a,Huang19b,BaiX23}.
In general, the gradient flow equation can be obtained from the thermodynamic variation of the free energy of the system
and is constrained by the dynamic dissipation path.
Typically, the total free energy $\mathcal{E}[\phi]$ arising from practical problems consists of
a quadratic term $\mathcal{E}_{_\mathcal{L}}[\phi]$ and a bounded, strong nonlinear term $\mathcal{E}_{_\mathcal{N}}[\phi]$, namely
\begin{equation}
\mathcal{E}[\phi]=\mathcal{E}_{_\mathcal{L}}[\phi]+\mathcal{E}_{_\mathcal{N}}[\phi],\quad
\text{where}~~\mathcal{E}_{_\mathcal{L}}[\phi]:=\frac{1}{2}(\phi,\mathcal{L}\phi),\label{eq1}
\end{equation}
with $\mathcal{L}$ is a linear, symmetric and non-negative operator, $(\cdot,\cdot)$ denotes the usual $L^{2}(\Omega)$-inner product.
Then, by performing the first variation of the free energy~\eqref{eq1} and applying the gradient flow idea,
we can arrive at the following dynamic equations
\begin{equation}
\left\{\begin{array}{l}
\frac{\partial\phi}{\partial t}=-\mathcal{G}\mu,\\[1mm]
\mu=\mathcal{L}\phi+H(\phi).
\end{array}\right.\label{eq2}
\end{equation}
Here, $\mathcal{G}$ is a positive operator related to the dynamic path, $\mu:=\delta \mathcal{E}/\delta\phi$ denotes the chemical potential
and $H(\phi):=\delta \mathcal{E}_{_\mathcal{N}}/\delta\phi$ represents the variational derivative of $\mathcal{E}_{_\mathcal{N}}[\phi]$.
Without loss of generality, Eq.~\eqref{eq2} is usually assumed to be equipped with suitable initial-boundary conditions,
making the problem complete. By using integration-by-parts with the boundary terms,
one can prove that the gradient flow of equation~\eqref{eq2} follows the important property of energy dissipation, i.e.,
\begin{equation}
\frac{\mathrm{d}\mathcal{E}[\phi]}{\mathrm{d}t}=\left(\frac{\delta \mathcal{E}}{\delta\phi},\frac{\partial\phi}{\partial t}\right)
=-(\mu,\mathcal{G}\mu)\le 0.\label{eq3}
\end{equation}

Due to the high nonlinearity of the term $H(\phi)$ in \eqref{eq2}, seeking its analytical solutions often becomes intractable.
Consequently, numerical methods with high accuracy are usually desirable in practical applications.
Meanwhile, in order to preserve the physical structure of the original, continuous system and conserve computational resources during long-term simulations,
researchers have been paying increasing attention on the development of structure-preserving numerical schemes, especially for the energy stability of the dissipative system~\eqref{eq3}.
To the best of our knowledge, three categories of unconditionally energy-stable methods have been developed for solving gradient flow systems \eqref{eq1}-\eqref{eq3}:
(i). the first type of methods discretize the original gradient flow equation \eqref{eq2} directly,
resulting in a strong nonlinear scheme that preserves the stability of the original energy,
such as the convex splitting method~\cite{Eyre98,Backofen19}, the exponential time differencing (ETD) method~\cite{Cox02,FuZ22}
and the discrete variational derivative (DVD) method~\cite{DuQ91,HuangJ20};
(ii). the second one transforms the original system into an equivalent system by introducing auxiliary variables which are dependent on the free energy,
then numerically discretizes this new system to achieve a linear scheme.
However, these methods can only ensure stability for a modified energy,
lacking an explicit connection between the auxiliary variables and their associated energy at the discrete level.
Examples include the invariant energy quadratization (IEQ) method~\cite{ChenC19,Huang19b},
the scalar auxiliary variable (SAV) method~\cite{ShenJ18,LiuZ20,HuangF20}
and the generalized positive auxiliary variable (gPAV) method~\cite{YangZ20,QianY20};
(iii). the third type of methods (e.g., the Lagrange multiplier-based SAV method~\cite{ChengQ20,Huang23,Huang23b},
the supplementary variable method (SVM)~\cite{GongY21,HongQ23} and the zero-factor method~\cite{LiuZ24}) are designed by adding
some new variables and incorporates the energy dissipation equation~\eqref{eq3} in numerical discretizations,
which usually result in a weakly nonlinear scheme, ensuring the stability of the original energy.

Among these methodologies, SAV-like approaches have garnered widespread acclaim in addressing a myriad of energy dissipation problems,
primarily attributed to their unconstrained free energy formulations and proficiency in devising high-order numerical schemes~\cite{ShenJ19,Huang22,Huang23}.
In terms of implementation paths, they are mainly divided into two categories: the nonlinear energy-based approach and the Lagrange multiplier approach.
The nonlinear energy-based SAV approach offers solutions irrespective of the time step size,
but it can only guarantee the stability of a modified energy at the discrete level~\cite{ShenJ18}.
Conversely, the Lagrange multiplier SAV approach guarantees the stability of the original energy \cite{ChengQ20},
albeit solutions may become infeasible when the time step size is excessively ``large'' \cite{Antoine21}.
In an endeavor to bridge the gap where the auxiliary variable is not inherently tied to its corresponding energy at the discrete level,
a relaxation technique has been introduced into the SAV method~\cite{JiangM22,Huang24}.
However, resolving the problem of potential non-existence of solutions in the Lagrange multiplier type case has not yet been addressed.

In this paper, we aim to develop a novel SAV framework that couples the two approaches described above and inherits their respective advantages
by incorporating weighting ideas.
In other words, the proposed method can have solutions for any time step size,
and establish an explicit connection between auxiliary variables and their associated energies at the discrete level,
maintaining the inherent advantages of the previous SAV approach.
In summary, the main contributions of this paper are as follows:
\begin{itemize}
  \item We unify the nonlinear energy-based SAV and Lagrange multiplier-based SAV approaches into a single framework for the first time;
  \item Our novel SAV framework establishes a direct connection between auxiliary variables and their associated energies at the discrete level,
        with which a discrete energy closer to the original one can be stabilized;
  \item The proposed weighted SAV approach is solvable for any time step size with sufficiently large weight coefficients.
\end{itemize}

The remainder of this paper is oganized as follows.
In Section~\ref{sec2}, we review the main concepts of the nonlinear energy-based SAV and Lagrange multiplier-based SAV approaches.
In Section~\ref{sec3}, a weighted SAV method is proposed by coupling the two types of SAV approaches, based on which the first-order and second-order numerical schemes are developed.
Several theorems are also presented to prove the existence of numerical solutions as well as energy stability. In Section~\ref{sec4}, we numerically verify the accuracy, efficiency and energy stability of our proposed schemes. Finally, some concluding remarks are given in Section~\ref{sec5}.

\section{Preliminaries}
\label{sec2}

In this section, we mainly review two types of SAV approaches widely used to solve gradient flow systems.

\subsection{Nonlinear energy-based SAV approach}

The key idea of nonlinear energy-based SAV approach~\cite{ShenJ18,LiuZ20,HuangF20} is to first introduce a scalar auxiliary variable $r:=r(t)$ that depends on the bounded nonlinear term $\mathcal{E}_{_\mathcal{N}}[\phi]$ in the free energy,
where the scalar auxiliary variable can be chosen from the power-type (e.g., the usual radical and/or linear styles),
exponential-type and logarithmic-type functions, and so on (see \cite{ChengQ21,Huang24} for details).
In the following, we review the main steps of the SAV approach by using the most popular radical-type scalar auxiliary variable~\cite{ShenJ18},
which is defined by
\begin{equation}
r(t)=\sqrt{\mathcal{E}_{_\mathcal{N}}[\phi]+C}.\label{eqa1}
\end{equation}
Here, the constant $C$ is used to ensure that the radicand is positive.

With the scalar auxiliary variable $r(t)$ defined in \eqref{eqa1}, the gradient flow equation~\eqref{eq2} can be rewritten as
\begin{equation}
\left\{\begin{array}{l}
\frac{\partial\phi}{\partial t}=-\mathcal{G}\mu,\\[1mm]
\mu=\mathcal{L}\phi+\frac{r}{\sqrt{\mathcal{E}_{_\mathcal{N}}[\phi]+C}}H(\phi),\\[2mm]
\frac{\mathrm{d}(r^{2}-C)}{\mathrm{d}t}=\frac{r}{\sqrt{\mathcal{E}_{_\mathcal{N}}[\phi]+C}}(H(\phi),\frac{\partial\phi}{\partial t}).
\end{array}\right.\label{eqa2}
\end{equation}
Next, taking the $L^{2}(\Omega)$-inner product of the first and second equations of \eqref{eqa2} with $\mu$ and $\phi_{t}$, respectively,
and combining the third equation leads to the following energy dissipation law
\begin{equation}
\frac{\mathrm{d}\widetilde{\mathcal{E}}[\phi,r]}{\mathrm{d}t}=-(\mu,\mathcal{G}\mu)\le 0,\quad\text{where}\quad
\widetilde{\mathcal{E}}[\phi,r]=\frac{1}{2}(\phi,\mathcal{L}\phi)+r^{2}-C.\label{eqa3}
\end{equation}
At the continuous level, the above energy $\widetilde{\mathcal{E}}[\phi,r]$, gradient flow equation~\eqref{eqa2} and energy dissipation law~\eqref{eqa3} dependent on the scalar auxiliary variable $r$ defined by \eqref{eqa1} are essentially the same as the original version,
which is also an important feature of SAV-like approaches.

In order to obtain the discrete equation associated with \eqref{eqa2} and satisfy the energy dissipation law,
we can use the special linear multi-step (e.g., backward Euler, Crank--Nicolson and second-order backward differentiation formulas)
or Runge--Kutta formulas to discretize the time variable~\cite{Huang23,Akrivis19,GongY20C,Feng21}.
To this end, we take the often-used backward Euler (BE) formula as an example to discretize the gradient flow equation
and obtain the following linear scheme:
\begin{equation}
\left\{\begin{array}{l}
\frac{\phi^{n+1}-\phi^{n}}{\tau }=-\mathcal{G}\mu^{n+1},\\[1mm]
\mu^{n+1}=\mathcal{L}\phi^{n+1}+\frac{r^{n+1}}{\sqrt{\mathcal{E}_{_\mathcal{N}}[\phi^{n}]+C}}H(\phi^{n}),\\[2mm]
2r^{n+1}\cdot\frac{r^{n+1}-r^{n}}{\tau }=\frac{r^{n+1}}{\sqrt{\mathcal{E}_{_\mathcal{N}}[\phi^{n}]+C}}
  \left(H(\phi^{n}),\frac{\phi^{n+1}-\phi^{n}}{\tau }\right),
\end{array}\right.\label{eqa4}
\end{equation}
where $\tau >0$ denotes the time step size and $\chi^{n} (\chi=\phi,\mu, r)$ represents the numerical approximation of $\chi(t)$
at time $t=t^{n}:=n\tau$ for $n=0,1,\ldots,N$ with $T=N\tau$.

For SAV-like approaches, the original variable $\phi, \mu$ and auxiliary variable $r$ are weakly coupled in the discrete equation,
which can be decoupled as follows.
Substituting the second equation into the first equation in \eqref{eqa4} and eliminating $\mu^{n+1}$, we obtain
\begin{equation}
(\mathcal{I}+\tau\mathcal{G}\mathcal{L})\phi^{n+1}=\phi^{n}-\frac{\tau\mathcal{G}H(\phi^{n})}{\sqrt{\mathcal{E}_{_\mathcal{N}}[\phi^{n}]+C}}r^{n+1},\label{eqa5}
\end{equation}
where $\mathcal{I}$ represents the identity operator.
Next, we let
\begin{equation}
\phi^{n+1}=\widetilde{p}^{n+1}+r^{n+1}\widetilde{q}^{n+1},\label{eqa6}
\end{equation}
according to the principle of linear superposition, $\widetilde{p}^{n+1}$ and $\widetilde{q}^{n+1}$ are the solutions of
\begin{equation}
(\mathcal{I}+\tau\mathcal{G}\mathcal{L})\widetilde{p}^{n+1}=\phi^{n}\quad\text{and}\quad
(\mathcal{I}+\tau\mathcal{G}\mathcal{L})\widetilde{q}^{n+1}=-\frac{\tau\mathcal{G}H(\phi^{n})}{\sqrt{\mathcal{E}_{_\mathcal{N}}[\phi^{n}]+C}},\label{eqa7}
\end{equation}
respectively.
Since $(\mathcal{I}+\tau\mathcal{G}\mathcal{L})$ is a positive operator and therefore invertible,
unique solutions exist for both of the above equations.
Once $\widetilde{p}^{n+1}$ and $\widetilde{q}^{n+1}$ are obtained, we can substitute \eqref{eqa6} into the last equation in \eqref{eqa4} to obtain the explicit result for $r^{n+1}$, i.e.,
\begin{equation}
r^{n+1}=\frac{2r^{n}\sqrt{\mathcal{E}_{_\mathcal{N}}[\phi^{n}]+C}+(H(\phi^{n}),\widetilde{p}^{n+1}-\phi^{n})}
{2\sqrt{\mathcal{E}_{_\mathcal{N}}[\phi^{n}]+C}-(H(\phi^{n}),\widetilde{q}^{n+1})},\label{eqa8}
\end{equation}
here we are only concerned with non-trivial solutions. 
Finally, the updated $\phi^{n+1}$ can be obtained by substituting $\widetilde{p}^{n+1}$, $\widetilde{q}^{n+1}$ and $r^{n+1}$ into \eqref{eqa6}.
As we can see, the primary computational effort of the aforementioned procedure lies in solving two linear elliptic partial differential equations with constant coefficients (i.e., Eq.~\eqref{eqa7}) at each time step.

Although the nonlinear energy-based SAV approach has many advantages such as simplicity, effectiveness and well-posedness,
it is somewhat regrettable that the scheme \eqref{eqa4} can only guarantee the energy dissipation for a modified energy:

\begin{thm}\label{thm1} (cf.~\cite{ShenJ18,ShenJ19})
The nonlinear energy-based SAV-BE scheme \eqref{eqa4} is unconditionally energy stable in the sense that
\begin{equation}
\frac{\widetilde{\mathcal{E}}[\phi^{n+1},r^{n+1}]-\widetilde{\mathcal{E}}[\phi^{n},r^{n}]}{\tau}\le -(\mu^{n+1},\mathcal{G}\mu^{n+1})\le 0,
\quad\text{where}\quad
\widetilde{\mathcal{E}}[\phi^{n},r^{n}]=\frac{1}{2}(\phi^{n},\mathcal{L}\phi^{n})+(r^{n})^{2}-C.
\label{tha1}
\end{equation}
\end{thm}

In discrete terms, the modified energy $\widetilde{\mathcal{E}}[\phi^{n},r^{n}]$ is merely a first-order approximation of
the original energy $\mathcal{E}[\phi^{n}]$, which means the above result does not ensure the dissipation of the original energy.
Actually, as demonstrated in previous numerical examples~\cite{ChenC19,JiangM22}, when the time step size $\tau$ is ``large'',
the original energy may become unstable even though the modified energy remains stable,
and there may be a large deviation between the modified energy and the original energy.
The main reason for this issue is that there is no explicit relationship between auxiliary variable $r^{n+1}$
and their associated energy $\sqrt{\mathcal{E}_{_\mathcal{N}}[\phi^{n+1}]+C}$ at the discrete level.
To solve this problem, a Lagrange multiplier strategy was introduced into the SAV framework to enforce dissipation of the original energy~\cite{ChengQ20}, which is described in the next subsection.

\subsection{Lagrange multiplier-based SAV approach}

By introducing a scalar Lagrange multiplier $\eta(t)$,
the gradient flow equation~\eqref{eq2} can be updated as~\cite{ChengQ20}
\begin{equation}
\left\{\begin{array}{l}
\frac{\partial\phi}{\partial t}=-\mathcal{G}\mu,\\[1mm]
\mu=\mathcal{L}\phi+\eta(t) H(\phi),\\[1mm]
\frac{\mathrm{d}\mathcal{E}_{_\mathcal{N}}[\phi]}{\mathrm{d}t}=\eta(t) (H(\phi),\frac{\partial\phi}{\partial t}).
\end{array}\right.\label{eqb1}
\end{equation}
If the initial value of the Lagrange multiplier is set as $\eta(0)=1$, the updated system above can be made equivalent to the original system, thus adhering to the original energy dissipation law \eqref{eq3}.

The continuous equation~\eqref{eqb1} can also be discretized by the backward Euler scheme in time:
\begin{equation}
\left\{\begin{array}{l}
\frac{\phi^{n+1}-\phi^{n}}{\tau }=-\mathcal{G}\mu^{n+1},\\[1mm]
\mu^{n+1}=\mathcal{L}\phi^{n+1}+\eta^{n+1}H(\phi^{n}),\\[1mm]
\frac{\mathcal{E}_{_\mathcal{N}}[\phi^{n+1}]-\mathcal{E}_{_\mathcal{N}}[\phi^{n}]}{\tau }
   =\eta^{n+1}\left(H(\phi^{n}),\frac{\phi^{n+1}-\phi^{n}}{\tau }\right),
\end{array}\right.\label{eqb2}
\end{equation}
where $\tau >0$ denotes the time step size and $\chi^{n}$ $(\chi=\phi,\mu,\eta)$ represents the numerical approximation of $\chi(t^{n})$.
Next, by eliminating $\mu^{n+1}$ in the first two equations in ~\eqref{eqb2}, we can obtain the following equation:
\begin{equation}
(\mathcal{I}+\tau\mathcal{G}\mathcal{L})\phi^{n+1}=\phi^{n}-\tau\mathcal{G}H(\phi^{n})\eta^{n+1},\label{eqb3}
\end{equation}
where $\mathcal{I}$ is the identity operator.
Same as before, we can decompose $\phi^{n+1}$ into the following two parts
\begin{equation}
\phi^{n+1}=\widehat{p}^{n+1}+\eta^{n+1}\widehat{q}^{n+1},\label{eqb4}
\end{equation}
where $\widehat{p}^{n+1}$ and $\widehat{q}^{n+1}$ are the solutions of
\begin{equation}
(\mathcal{I}+\tau\mathcal{G}\mathcal{L})\widehat{p}^{n+1}=\phi^{n}\quad\text{and}\quad
(\mathcal{I}+\tau\mathcal{G}\mathcal{L})\widehat{q}^{n+1}=-\tau\mathcal{G}H(\phi^{n}),\label{eqb5}
\end{equation}
which are uniquely solvable equations.
Finally, by substituting \eqref{eqb4} into the third equation in \eqref{eqb2}, we have
\begin{equation}
\mathcal{E}_{_\mathcal{N}}[\widehat{p}^{n+1}+\eta^{n+1}\widehat{q}^{n+1}]
-(\eta^{n+1})^{2}(H(\phi^{n}),\widehat{q}^{n+1})
-\eta^{n+1}(H(\phi^{n}),\widehat{p}^{n+1}-\phi^{n})
-\mathcal{E}_{_\mathcal{N}}[\phi^{n}]=0,\label{eqb6}
\end{equation}
leading to a nonlinear scalar equation with respect to $\eta^{n+1}$. Its degree of nonlinearity mainly depends on $\mathcal{E}_{_\mathcal{N}}$.
In particular, since the Lagrange multiplier $\eta(t)\equiv 1$ at the continuous level,
we focus on the solution of Eq.~\eqref{eqb6} near $1$.
Recently, theoretical results~\cite{Onuma24} have proved that under certain mild conditions,
Eq.~\eqref{eqb6} has a solution near $1$ when the time step size $\tau$ is small enough.
However, if the time step size $\tau$ is ``large'',
the numerical results presented in~\cite{Antoine21} indicate that Eq.~\eqref{eqb6} may not converge to a solution near $1$,
or possibly have no solution at all.
Fortunately, when the scheme~\eqref{eqb2} does have a solution, its original energy remains stable.

\begin{thm}\label{thm2} (cf.~\cite{ChengQ20,Huang23})
If the Lagrange multiplier-based SAV-BE scheme \eqref{eqb2} has solutions,
then it satisfies the following original energy dissipation law
\begin{equation}
\frac{\mathcal{E}[\phi^{n+1}]-\mathcal{E}[\phi^{n}]}{\tau}\le -(\mu^{n+1},\mathcal{G}\mu^{n+1})\le 0,
\quad\text{where}\quad
\mathcal{E}[\phi^{n}]=\frac{1}{2}(\phi^{n},\mathcal{L}\phi^{n})+\mathcal{E}_{_\mathcal{N}}[\phi^{n}].
\label{thb1}
\end{equation}
\end{thm}

In summary, the main computational effort of SAV-like approaches, whether nonlinear energy-based type or Lagrange multiplier type,
only need to solve, twice, a decoupled linear, elliptic partial differential equations with constant coefficients at each time step.
On the other hand, they have their own disadvantages, i.e.,  the nonlinear energy-based SAV approach may result in a significant deviation between the auxiliary variable and its associated energy at the discrete level,
while the Lagrange multiplier-based SAV approach may have no solution for a ``large'' time step size.
In the subsequent section, with the aid of the weighting idea,
we consolidate the aforementioned two types of SAV approaches into a unified framework,
establishing an explicit relationship between auxiliary variables and their associated energies at the discrete level,
and ensuring that proposed schemes have solutions for arbitrary time step sizes.

\section{Weighted SAV approach}
\label{sec3}

In this section, we propose a novel method that utilizes the concept of weighting to unify the two distinct types of SAV approaches presented above into a cohesive framework,
which will be discussed in detail below.

Considering the scalar auxiliary variable $r(t)=\sqrt{\mathcal{E}_{_\mathcal{N}}[\phi]+C}$ introduced in \eqref{eqa1},
we would adopt a new reconstruction style of the gradient flow equation \eqref{eq2} as follows
\begin{equation}
\left\{\begin{array}{l}
\frac{\partial\phi}{\partial t}=-\mathcal{G}\mu,\\[1mm]
\mu=\mathcal{L}\phi+\frac{r}{\sqrt{\mathcal{E}_{_\mathcal{N}}[\phi]+C}}H(\phi),\\[2.5mm]
\lambda\frac{\mathrm{d}(r^{2}-C)}{\mathrm{d}t}+(1-\lambda)\frac{\mathrm{d}\mathcal{E}_{_\mathcal{N}}[\phi]}{\mathrm{d}t}
 =\frac{r}{\sqrt{\mathcal{E}_{_\mathcal{N}}[\phi]+C}}\big(H(\phi),\frac{\partial\phi}{\partial t}\big),
\end{array}\right.\label{eqc1}
\end{equation}
where $\lambda\in[0,1]$ is a weight coefficient.
It is clearly observed that, this weighted SAV approach is of nonlinear energy-based type when $\lambda=1$;
and when $\lambda=0$, \eqref{eqc1} is regarded as an equivalent formulation of Lagrange multiplier-based SAV, since we can view $r(t)/\sqrt{\mathcal{E}_{_\mathcal{N}}[\phi]+C}$ as an alternative form of $\eta(t)$.
In other words, the above reconstruction equation~\eqref{eqc1} includes the nonlinear energy-based SAV
and Lagrange multiplier-based SAV approaches introduced earlier, as well as their intermediate states.

Similarly, by taking the $L^{2}(\Omega)$-inner product of the first and second equations of \eqref{eqc1} with $\mu$ and $\phi_{t}$, respectively,
and combining the third equation, we can derive the following energy dissipation law,
\begin{equation}
\frac{\mathrm{d}\overline{\mathcal{E}}[\phi,r]}{\mathrm{d}t}=-(\mu,\mathcal{G}\mu)\le 0,\quad\text{where}\quad
\overline{\mathcal{E}}[\phi,r]=\frac{1}{2}(\phi,\mathcal{L}\phi)+\lambda(r^{2}-C)+(1-\lambda)\mathcal{E}_{_\mathcal{N}}[\phi].\label{eqc2}
\end{equation}
At the continuous level, when $r^{2}-C\equiv\mathcal{E}_{_\mathcal{N}}[\phi]$,
the gradient flow equation~\eqref{eqc1} and its associated energy dissipation law~\eqref{eqc2}
are essentially the same as the original style regardless the weight coefficient $\lambda$.
At the discrete level, as we will see,
a smaller weight coefficient $\lambda$ closer to zero would make the modified energy closer to the original energy.
On the other hand, a very small value of $\lambda$ may lead to the weighted SAV scheme reducing to the Lagrange multiplier-based case,
potentially resulting in no solution when the time step size $\tau$ is sufficiently large.
Consequently, the key to the success of this weighted SAV method lies in determining an optimal weight $\lambda_{opt}$ that is as small as possible, yet still guarantees the existence of a solution, particularly in scenarios involving large time steps.
In the following, we will first introduce the discretization of~\eqref{eqc1}
and then develop some practical methods in searching such weight coefficient $\lambda$.

\subsection{First-order weighted SAV scheme}

To begin with, we will utilize the backward Euler (BE) scheme to discretize the above continuous equation~\eqref{eqc1}:
\begin{equation}
\left\{\begin{array}{l}
\frac{\phi^{n+1}-\phi^{n}}{\tau }=-\mathcal{G}\mu^{n+1},\\[1mm]
\mu^{n+1}=\mathcal{L}\phi^{n+1}+\frac{r^{n+1}}{\sqrt{\mathcal{E}_{_\mathcal{N}}[\phi^{n}]+C}}H(\phi^{n}),\\[2mm]
\lambda\cdot2r^{n+1}\cdot\frac{r^{n+1}-r^{n}}{\tau }
 +(1-\lambda)\frac{\mathcal{E}_{_\mathcal{N}}[\phi^{n+1}]-\mathcal{E}_{_\mathcal{N}}[\phi^{n}]}{\tau}
 =\frac{r^{n+1}}{\sqrt{\mathcal{E}_{_\mathcal{N}}[\phi^{n}]+C}}\left(H(\phi^{n}),\frac{\phi^{n+1}-\phi^{n}}{\tau }\right),
\end{array}\right.\label{eqc3}
\end{equation}
where $\tau >0$ represents the time step size.
Different from the case of nonlinear energy-based type ($\lambda=1$),
in which the auxiliary variable $r$ only has an explicit relationship with $\mathcal{E}_{_\mathcal{N}}[\phi]$ at the continuous level,
the auxiliary variable $r^{n+1}$ also has a direct connection with its associated energy $\mathcal{E}_{_\mathcal{N}}[\phi^{n+1}]$ at the discrete level, as described by the last equation of \eqref{eqc3}.
In the meanwhile, since the first two equations in \eqref{eqc3} and \eqref{eqa4} are exactly the same,
we can take a similar decoupling strategy for $\phi^{n+1}$ and $r^{n+1}$ as follows.
Firstly, we eliminate $\mu^{n+1}$ and obtain
\begin{equation}
(\mathcal{I}+\tau\mathcal{G}\mathcal{L})\phi^{n+1}=\phi^{n}-\frac{\tau\mathcal{G}H(\phi^{n})}{\sqrt{\mathcal{E}_{_\mathcal{N}}[\phi^{n}]+C}}r^{n+1},
\label{eqc4}
\end{equation}
where $\mathcal{I}$ is defined as the identity operator.
Next, we decompose $\phi^{n+1}$ into the following form
\begin{equation}
\phi^{n+1}=p^{n+1}+r^{n+1}q^{n+1},\label{eqc5}
\end{equation}
where $p^{n+1}$ and $q^{n+1}$ are the solutions of
\begin{equation}
(\mathcal{I}+\tau\mathcal{G}\mathcal{L})p^{n+1}=\phi^{n}\quad\text{and}\quad
(\mathcal{I}+\tau\mathcal{G}\mathcal{L})q^{n+1}=-\frac{\tau\mathcal{G}H(\phi^{n})}{\sqrt{\mathcal{E}_{_\mathcal{N}}[\phi^{n}]+C}},\label{eqc6}
\end{equation}
respectively. The invertibility of $(\mathcal{I}+\tau\mathcal{G}\mathcal{L})$ for any $\tau>0$ makes the above two equations uniquely solvable.
Finally, by substituting $p^{n+1}$ and $q^{n+1}$ into \eqref{eqc5} and putting it into the last equation in \eqref{eqc3}, we have the following nonlinear scalar equation with respect to $r^{n+1}$, for a given $\lambda\in[0,1]$:
\begin{equation}
2\lambda r^{n+1}(r^{n+1}-r^{n})+(1-\lambda)\left(\mathcal{E}_{_\mathcal{N}}[p^{n+1}+ r^{n+1}q^{n+1}]-\mathcal{E}_{_\mathcal{N}}[\phi^{n}]\right)
-\frac{r^{n+1}(H(\phi^{n}),p^{n+1}+ r^{n+1}q^{n+1}-\phi^{n})}{\sqrt{\mathcal{E}_{_\mathcal{N}}[\phi^{n}]+C}}=0,\label{eqc7}
\end{equation}

Now we come to the crucial step of selecting the weight coefficient $\lambda$.
Before we introduce the strategy in this selection process,
we would first present the following results on the existence of solutions and energy stability of the weighted SAV-BE scheme~\eqref{eqc3}.

\begin{thm}\label{thmf1}
For any given time step size $\tau>0$, the weighted SAV-BE scheme \eqref{eqc3} ensures the existence of solutions when the weight coefficients are sufficiently large. More specifically, there exists a constant $\underline{\lambda}\in(0,1)$ such that the scheme always has solutions for all $\lambda\in(\underline{\lambda},1]$.
\end{thm}

\begin{proof}
For the weighted SAV-BE scheme \eqref{eqc3}, $\lambda=1$ and $\lambda=0$ correspond to nonlinear energy-based type and Lagrange multiplier type cases, respectively, where the former has solutions for any time step size $\tau>0$,
while the latter has solutions when $\tau$ is sufficiently small~\cite{Onuma24}.
In this theorem, we focus on the case of $0<\lambda<1$.

Since $\mathcal{G}$ and $\mathcal{L}$ are positive and non-negative operators, respectively,
it follows that $(\mathcal{I}+\tau\mathcal{G}\mathcal{L})$ is a positive operator for any $\tau>0$ and thus invertible.
This means that both equations of $p^{n+1}$ and $q^{n+1}$ in \eqref{eqc6} have unique solutions,
leaving us only to consider the solvability of nonlinear scalar equations \eqref{eqc7} for $r^{n+1}$.

For the sake of simplicity, based on Eq.~\eqref{eqc7}, we introduce the following two expressions
\begin{align}
f(x;\lambda):\;=&\, 2\lambda x(x-r^{n})
-\frac{x}{\sqrt{\mathcal{E}_{_\mathcal{N}}[\phi^{n}]+C}}(H(\phi^{n}),p^{n+1}+ q^{n+1}x-\phi^{n})\nonumber\\
=&\left[2\lambda-\frac{(H(\phi^{n}),q^{n+1})}{\sqrt{\mathcal{E}_{_\mathcal{N}}[\phi^{n}]+C}}\right]x^{2}
-\left[2 r^{n}\lambda+\frac{(H(\phi^{n}),p^{n+1}-\phi^{n})}{\sqrt{\mathcal{E}_{_\mathcal{N}}[\phi^{n}]+C}}\right]x\nonumber\\
=&:f_{2}(\lambda)x^{2}-f_{1}(\lambda)x,\label{thms1}
\end{align}
and
\begin{align}
g(x;\lambda):= (1-\lambda)\left(\mathcal{E}_{_\mathcal{N}}[p^{n+1}+ q^{n+1}x]-\mathcal{E}_{_\mathcal{N}}[\phi^{n}]\right).\label{thms2}
\end{align}

Firstly, we prove that the quadratic term coefficient $f_{2}(\lambda)$ of $f(x;\lambda)$ is always positive for any $\tau>0$.
The second equation in \eqref{eqc6} yields
\begin{equation}
\left(\mathcal{G}^{-1}+\tau\mathcal{L}\right)q^{n+1}=-\frac{\tau H(\phi^{n})}{\sqrt{\mathcal{E}_{_\mathcal{N}}[\phi^{n}]+C}},\label{thms3}
\end{equation}
then, taking the $L^{2}(\Omega)$-inner product of both ends of the above equation with $q^{n+1}$, we obtain
\begin{equation}
-\frac{1}{\sqrt{\mathcal{E}_{_\mathcal{N}}[\phi^{n}]+C}}(H(\phi^{n}),q^{n+1})
=\frac{1}{\tau}\left((\mathcal{G}^{-1}+\tau\mathcal{L})q^{n+1},q^{n+1}\right)\ge 0,\label{thms4}
\end{equation}
where the last inequality is due to the fact that $\mathcal{G}$ and $\mathcal{L}$ are positive and non-negative operators, respectively.
Since $\lambda>0$, it follows that $f_{2}(\lambda)>0$.

Secondly, we discuss the sign of the linear term coefficient $f_{1}(\lambda)$ of $f(x;\lambda)$.
If $f_{1}(\lambda)$ has a zero-point $\lambda_{0}$ within the interval $(0,1)$,
then according to the continuity of $f_{1}(\lambda)$, it can be inferred that $f_{1}(\lambda)\neq 0$ when $\lambda>\lambda_{0}$.
Conversely, if $f_{1}(\lambda)$ has no zero-point within the interval $(0,1)$,
then for any $\lambda\in(0,1)$ there is $f_{1}(\lambda)\neq 0$.
All in all, there exists $0\leq\lambda_{0}<1$ such that $f_{1}(\lambda)\neq 0$ for any $\lambda_{0}<\lambda<1$.
Consequently, we may conclude that $f(x;\lambda)$ must have a negative minimum $f_{\min}$ for any $\lambda\in (\lambda_{0},1)$,
since the quadratic term coefficient $f_{2}(\lambda)$ of $f(x;\lambda)$ is always positive
and the linear term coefficient $f_{1}(\lambda)$ is non-zero when $\lambda>\lambda_{0}$.
Moreover, considering that $\mathcal{E}_{_\mathcal{N}}$ is bounded,
there must exist a positive constant $M$ such that $|\mathcal{E}_{_\mathcal{N}}|\le M$.
To ensure that Eq.~\eqref{eqc7} has solutions, it is sufficient to satisfy
\begin{align}
[f(x;\lambda)+g(x;\lambda)]_{\min}
\le [f(x;\lambda)]_{\min}+[g(x;\lambda)]_{\max}\le f_{\min}+2M(1-\lambda) \le 0,\label{thms5}
\end{align}
that is, $\lambda\ge 1+\frac{f_{\min}}{2M}$.
As a result, by choosing $\underline{\lambda}=\max\{\lambda_{0},1+\frac{f_{\min}}{2M}\}$,
then for any $\lambda\in(\underline{\lambda},1]$,
the nonlinear scalar equation \eqref{eqc7} always has solutions. The proof is completed.
\end{proof}

The above results illustrate that within the framework of SAV approach,
we can employ the nonlinear energy-based type to refine the Lagrange multiplier type, thereby ensuring that the proposed numerical scheme~\eqref{eqc3} has solutions for any time step size $\tau>0$.
On the other hand, as indicated by the following theorem, weighting with the Lagrange multiplier type may help to adjust the discrete energy to be closer to the original energy.
Under certain conditions, moreover, it can even ensure the stability of the original energy.

\begin{thm}\label{thmf2}
The weighted SAV-BE scheme \eqref{eqc3} unconditionally satisfies the following energy dissipation law:
\begin{equation}
\frac{\overline{\mathcal{E}}[\phi^{n+1},r^{n+1}]-\overline{\mathcal{E}}[\phi^{n},r^{n}]}{\tau}\le -(\mu^{n+1},\mathcal{G}\mu^{n+1})\le 0,\label{thc1}
\end{equation}
where
\begin{equation}
\overline{\mathcal{E}}[\phi^{n},r^{n}]=\frac{1}{2}(\phi^{n},\mathcal{L}\phi^{n})
  +\lambda((r^{n})^{2}-C)+(1-\lambda)\mathcal{E}_{_\mathcal{N}}[\phi^{n}].\label{thc2}
\end{equation}
Furthermore, if $2r^{n+1}(r^{n+1}-r^{n})\ge \mathcal{E}_{_\mathcal{N}}[\phi^{n+1}]-\mathcal{E}_{_\mathcal{N}}[\phi^{n}]$,
the proposed scheme also complies with the following original energy dissipation law:
\begin{equation}
\frac{\mathcal{E}[\phi^{n+1}]-\mathcal{E}[\phi^{n}]}{\tau}\le -(\mu^{n+1},\mathcal{G}\mu^{n+1})\le 0,
\quad\text{where}\quad
\mathcal{E}[\phi^{n}]=\frac{1}{2}(\phi^{n},\mathcal{L}\phi^{n})+\mathcal{E}_{_\mathcal{N}}[\phi^{n}].\label{thc3}
\end{equation}
\end{thm}

\begin{proof}
Taking the $L^{2}(\Omega)$-inner products of \eqref{eqc3}$_{1}$ and \eqref{eqc3}$_{2}$ with $\mu^{n+1}$ and $\frac{\phi^{n+1}-\phi^{n}}{\tau}$, respectively,
we obtain
\begin{align}
&\left(\frac{\phi^{n+1}-\phi^{n}}{\tau },\mu^{n+1}\right)=-(\mathcal{G}\mu^{n+1},\mu^{n+1}),\label{thc4}\\
&\left(\mu^{n+1},\frac{\phi^{n+1}-\phi^{n}}{\tau }\right)
=\left(\mathcal{L}\phi^{n+1},\frac{\phi^{n+1}-\phi^{n}}{\tau }\right)
  +\frac{r^{n+1}}{\sqrt{\mathcal{E}_{_\mathcal{N}}[\phi^{n}]+C}}\left(H(\phi^{n}),\frac{\phi^{n+1}-\phi^{n}}{\tau }\right),\label{thc5}
\end{align}
then, combination \eqref{eqc3}$_{3}$ results in
\begin{align}
-\tau(\mathcal{G}\mu^{n+1},\mu^{n+1})
=&\left(\mathcal{L}\phi^{n+1},\phi^{n+1}-\phi^{n}\right)
  +\frac{r^{n+1}}{\sqrt{\mathcal{E}_{_\mathcal{N}}[\phi^{n}]+C}}\left(H(\phi^{n}),\phi^{n+1}-\phi^{n}\right)\nonumber\\
=&\left(\mathcal{L}\phi^{n+1},\phi^{n+1}-\phi^{n}\right)
  +\lambda\cdot2r^{n+1}(r^{n+1}-r^{n})+(1-\lambda)(\mathcal{E}_{_\mathcal{N}}[\phi^{n+1}]-\mathcal{E}_{_\mathcal{N}}[\phi^{n}])\nonumber\\
=&\,\frac{1}{2}\left(\mathcal{L}\phi^{n+1},\phi^{n+1}\right)-\frac{1}{2}\left(\mathcal{L}\phi^{n},\phi^{n}\right)
  +\frac{1}{2}\left(\mathcal{L}(\phi^{n+1}-\phi^{n}),\phi^{n+1}-\phi^{n}\right)\nonumber\\[1mm]
  &+\lambda\left[(r^{n+1})^{2}-(r^{n})^{2}+(r^{n+1}-r^{n})^{2}\right]
 +(1-\lambda)(\mathcal{E}_{_\mathcal{N}}[\phi^{n+1}]-\mathcal{E}_{_\mathcal{N}}[\phi^{n}]),\label{thc6}
\end{align}
the rearrangement of the above expression yields the energy dissipation law~\eqref{thc1} with \eqref{thc2}.

Furthermore, if $2r^{n+1}(r^{n+1}-r^{n})\ge \mathcal{E}_{_\mathcal{N}}[\phi^{n+1}]-\mathcal{E}_{_\mathcal{N}}[\phi^{n}]$,
it can be further obtained from \eqref{thc6} that
\begin{align}
-\tau(\mathcal{G}\mu^{n+1},\mu^{n+1})
=&\left(\mathcal{L}\phi^{n+1},\phi^{n+1}-\phi^{n}\right)
  +\lambda\cdot2r^{n+1}(r^{n+1}-r^{n})+(1-\lambda)(\mathcal{E}_{_\mathcal{N}}[\phi^{n+1}]-\mathcal{E}_{_\mathcal{N}}[\phi^{n}])\nonumber\\
\ge&\left(\mathcal{L}\phi^{n+1},\phi^{n+1}-\phi^{n}\right)
  +\lambda(\mathcal{E}_{_\mathcal{N}}[\phi^{n+1}]-\mathcal{E}_{_\mathcal{N}}[\phi^{n}])
  +(1-\lambda)(\mathcal{E}_{_\mathcal{N}}[\phi^{n+1}]-\mathcal{E}_{_\mathcal{N}}[\phi^{n}])\nonumber\\
\ge&\,\frac{1}{2}\left(\mathcal{L}\phi^{n+1},\phi^{n+1}\right)-\frac{1}{2}\left(\mathcal{L}\phi^{n},\phi^{n}\right)
  +\mathcal{E}_{_\mathcal{N}}[\phi^{n+1}]-\mathcal{E}_{_\mathcal{N}}[\phi^{n}],\label{thc7}
\end{align}
which implies the original energy dissipation law~\eqref{thc3}.
\end{proof}

As shown in Eqs.~\eqref{thc1}-\eqref{thc2} of the previous theorem,
a smaller value of $\lambda$ leads to a dissipation of discrete energy closer to $\mathcal{E}[\phi^{n}]$.
However, if $\lambda$ is too small, it may not ensure the existence of a solution as stipulated by Theorem \ref{thmf1}.
To achieve such an optimal $\lambda$ in practical simulations, we need the following critical observation, which reveals the impact of the quantitative relationship between $2r^{n+1}(r^{n+1}-r^{n})$ and $\mathcal{E}_{_\mathcal{N}}[\phi^{n+1}]-\mathcal{E}_{_\mathcal{N}}[\phi^{n}]$ on the existence of solutions for \eqref{eqc3}.

\begin{thm}\label{thmf3}
If the weight coefficient $\lambda=\lambda_{\star}\in[0,1]$ yields a solution $(\phi_{\star}^{n+1},\mu_{\star}^{n+1},r_{\star}^{n+1})$ for the weighted SAV-BE scheme \eqref{eqc3}, then

(i) when $2r_{\star}^{n+1}(r_{\star}^{n+1}-r^{n})\ge \mathcal{E}_{_\mathcal{N}}[\phi_{\star}^{n+1}]-\mathcal{E}_{_\mathcal{N}}[\phi^{n}]$,
the scheme \eqref{eqc3} always has solutions for any $\lambda\in[0,\lambda_{\star}]$;

(ii) when $2r_{\star}^{n+1}(r_{\star}^{n+1}-r^{n})<\mathcal{E}_{_\mathcal{N}}[\phi_{\star}^{n+1}]-\mathcal{E}_{_\mathcal{N}}[\phi^{n}]$,
the scheme \eqref{eqc3} always has solutions for any $\lambda\in[\lambda_{\star},1]$.
\end{thm}

\begin{proof}
Based on the nonlinear scalar equation \eqref{eqc7}, we introduce the following bivariate function of $(r,\lambda)$:
\begin{align}
h(r,\lambda):=&\,2\lambda r(r-r^{n})
+(1-\lambda)\left(\mathcal{E}_{_\mathcal{N}}[p^{n+1}+q^{n+1}r]-\mathcal{E}_{_\mathcal{N}}[\phi^{n}]\right)\nonumber\\
&\,-\frac{r}{\sqrt{\mathcal{E}_{_\mathcal{N}}[\phi^{n}]+C}}(H(\phi^{n}),p^{n+1}+q^{n+1}r-\phi^{n})\nonumber\\
=&\left[2 r(r-r^{n})-(\mathcal{E}_{_\mathcal{N}}[p^{n+1}+q^{n+1}r]-\mathcal{E}_{_\mathcal{N}}[\phi^{n}])\right]\lambda\nonumber\\
&\,+\left(\mathcal{E}_{_\mathcal{N}}[p^{n+1}+q^{n+1}r]-\mathcal{E}_{_\mathcal{N}}[\phi^{n}]\right)
-\frac{r}{\sqrt{\mathcal{E}_{_\mathcal{N}}[\phi^{n}]+C}}(H(\phi^{n}),p^{n+1}+q^{n+1}r-\phi^{n})\nonumber\\
=&\left[2\lambda-\frac{(H(\phi^{n}),q^{n+1})}{\sqrt{\mathcal{E}_{_\mathcal{N}}[\phi^{n}]+C}}\right]r^{2}
-\left[2 r^{n}\lambda+\frac{(H(\phi^{n}),p^{n+1}-\phi^{n})}{\sqrt{\mathcal{E}_{_\mathcal{N}}[\phi^{n}]+C}}\right]r\nonumber\\
&\,+(1-\lambda)\left(\mathcal{E}_{_\mathcal{N}}[p^{n+1}+q^{n+1}r]-\mathcal{E}_{_\mathcal{N}}[\phi^{n}]\right).\label{thmh1}
\end{align}
We first consider case (i). Noting that $p^{n+1}+q^{n+1}r_{\star}^{n+1}=\phi_{\star}^{n+1}$,
when $2r_{\star}^{n+1}(r_{\star}^{n+1}-r^{n})\ge \mathcal{E}_{_\mathcal{N}}[\phi_{\star}^{n+1}]-\mathcal{E}_{_\mathcal{N}}[\phi^{n}]$,
the univariate function $h(r_{\star}^{n+1},\lambda)$ is monotonically increasing with respect to $\lambda$.
Therefore, for any $\lambda_{1}\in[0,\lambda_{\star}]$, we have
\begin{equation}
h(r_{\star}^{n+1},\lambda_{1})\le h(r_{\star}^{n+1},\lambda_{\star})=0.\label{thmh2}
\end{equation}
From \eqref{thms4}, we know that $(H(\phi^{n}),q^{n+1})\le 0$.
Therefore, when $\lambda_{1}\in(0,\lambda_{\star}]$ or $\lambda_{1}=0$ and $(H(\phi^{n}),q^{n+1})<0$,
the coefficient of $r^{2}$ in $h(r,\lambda_{1})$ is always positive.
Considering that $\mathcal{E}_{_\mathcal{N}}$ is bounded, it follows that
\begin{equation}
\lim_{r\to -\infty} h(r,\lambda_{1})=+\infty
\quad\text{and}\quad
\lim_{r\to +\infty} h(r,\lambda_{1})=+\infty.\label{thmh3}
\end{equation}
According to the zero-point theorem,
we know that the equation $h(r,\lambda_{1})=0$ of $r$ has solutions
within the intervals $(-\infty,r_{\star}^{n+1}]$ and $[r_{\star}^{n+1},+\infty)$.

Next, we discuss the case where $\lambda_{1}=0$ and $(H(\phi^{n}),q^{n+1})=0$,
in which the coefficient of $r^{2}$ in $h(r,\lambda_{1})$ is always zero.
Under this premise, if $(H(\phi^{n}),p^{n+1}-\phi^{n})<0$,
then the coefficient of $r$ in $h(r,\lambda_{1})$ is always positive, resulting in $\lim_{r\to +\infty} h(r,\lambda_{1})=+\infty$.
This implies that the equation $h(r,\lambda_{1})=0$ of $r$ has solutions within $[r_{\star}^{n+1},+\infty)$.
Similarly, if $(H(\phi^{n}),p^{n+1}-\phi^{n})>0$, then the equation $h(r,\lambda_{1})=0$ of $r$ has solutions within $(-\infty,r_{\star}^{n+1}]$.
Finally, if $(H(\phi^{n}),p^{n+1}-\phi^{n})=0$, then $r_{\star}^{n+1}$ is the solution to the equation $h(r,\lambda_{1})=0$ of $r$.
This completes the proof for case (i).


Similarly, we can prove case (ii). Thus, the theorem is proved.
\end{proof}

The above results motivate us to follow the procedure outlined below to determine the optimal value of $\lambda_{opt}$:
at the $n$-th step, we first verify if Eq.~\eqref{eqc7} with $\lambda=0$ has a solution, and select $\lambda_{opt}=0$ if it does.
Otherwise, we must have
$2r_{\star}^{n+1}(r_{\star}^{n+1}-r^{n})<\mathcal{E}_{_\mathcal{N}}[\phi_{\star}^{n+1}]-\mathcal{E}_{_\mathcal{N}}[\phi^{n}]$
according to (ii) in Theorem \ref{thmf3}.
In this case, following Theorem \ref{thmf1}, we would evaluate the solution existence of Eq.~\eqref{eqc7} by assuming $\lambda$ to be some positive number $\lambda_{a}$ (e.g., 0.5).
In the case that \eqref{eqc7} can be solved with $\lambda_a$, we can further assume that $\lambda_{opt}\in [0,\lambda_a]$,
since we are seeking the minimum value of $\lambda$ and result (ii) in Theorem \ref{thmf3} indicates that \eqref{eqc7} always has solution for any $\lambda\in[\lambda_a,1]$.
On the other hand, based on the result (ii) in Theorem \ref{thmf3} again,
we may conclude that there is no $\lambda<\lambda_a$ that guarantees the existence of a solution for Eq.~\eqref{eqc7}.
In such a scenario, we would proceed by making a new guess of $\lambda_{opt}$ within the range $[\lambda_a,1]$.

\begin{algorithm}
\caption{The first-order weighted SAV-BE scheme \eqref{eqc3}.}\label{alg1}
Update the unknown variables $(\phi^{n+1},r^{n+1})$ by the known ones $(\phi^{n},r^{n})$.
\begin{enumerate}
  \item Solve $p^{n+1}$ and $q^{n+1}$ form \eqref{eqc6}.
  \item Set $\lambda =0$ and solve Eq.~\eqref{eqc7}.\\
  \textbf{if} Eq.~\eqref{eqc7} is solvable \textbf{then} \\
  $~~~~$ Output the final values of $r^{n+1}$ and $\lambda=0$.\\
  \textbf{else}  \\
  $~~~~$ Set $a=0, b=1$ and provide the tolerance $tol$.\\
  $~~~~$ \textbf{while} $b-a\ge tol$ \textbf{do}\\
  $~~~~~~~~~$ Set $\lambda=\frac{b+a}{2}$ and solve Eq.~\eqref{eqc7}.\\
  $~~~~~~~~~$ \textbf{if} Eq.~\eqref{eqc7} is solvable \textbf{then}\\
  $~~~~~~~~~~~~~$ Set $a=a$ and $b=\frac{b+a}{2}$.\\
  $~~~~~~~~~$ \textbf{else}\\
  $~~~~~~~~~~~~~$ Set $a=\frac{b+a}{2}$ and $b=b$.\\
  $~~~~~~~~~$ \textbf{end if}\\
  $~~~~$ \textbf{end while}\\
  $~~~~$ Output the final values of $r^{n+1}$ and $\lambda$.\\
  \textbf{end if}
  \item Update $\phi^{n+1}$ by \eqref{eqc5}.
\end{enumerate}
\end{algorithm}

By repeatedly applying this idea, we can formulate Algorithm \ref{alg1} for the weighted SAV-BE scheme \eqref{eqc3}.
In this algorithm, in which a straightforward bisection method is applied to determine the optimal value of $\lambda$.
Moreover, considering the nonlinearity of $\mathcal{E}_{_\mathcal{N}}[\phi]$
and the role of $r^{n+1}/\sqrt{\mathcal{E}_{_\mathcal{N}}[\phi^{n}]+C}$ akin to the Lagrange multiplier $\eta^{n+1}$,
we can employ Newton's method to solve the nonlinear scalar equation \eqref{eqc7}
with an initial guess of $\sqrt{\mathcal{E}_{_\mathcal{N}}[\phi^{n}]+C}$.
In summary, the weighted SAV-BE scheme \eqref{eqc3} needs to solve two elliptic partial differential equations with constant coefficients
and several nonlinear scalar equations at each time step, in which the main computational cost is consistent with the classical SAV scheme.

\begin{remark}\label{remF}
In fields such as materials science, the energy density of $\mathcal{E}_{_\mathcal{N}}[\phi]$ is often considered in the form of a fourth-order polynomial with a double-well structure.
This simplifies Eq.~\eqref{eqc7} to a quartic equation in one variable,
enabling us to easily apply the discriminant to determine the nature of its roots.
Moreover, to avoid the scenario where Eq.~\eqref{eqc7} may have no solution due to solver issues,
we will further analyze its monotonicity to reduce the risk of misjudgment.
\end{remark}

\begin{remark}\label{remOpt}
The core challenge of the aforementioned weighted SAV approach, as we have mentioned,
is how to determine the minimum value of the weighting coefficient $\lambda$ given a specified time step size.
An alternative approach to address this issue can be formulated as the following constrained optimization problem:
\begin{align}
&\lambda=\min\limits_{\overline{\lambda}\in[0,1],\,\phi^{n+1},\,\mu^{n+1},\,r^{n+1}}\overline{\lambda},\label{eqrema}\\[1mm]
&\text{s.t.}~
\left\{\begin{array}{l}
\frac{\phi^{n+1}-\phi^{n}}{\tau }=-\mathcal{G}\mu^{n+1},\\[1.5mm]
\mu^{n+1}=\mathcal{L}\phi^{n+1}+\frac{r^{n+1}}{\sqrt{\mathcal{E}_{_\mathcal{N}}[\phi^{n}]+C}}H(\phi^{n}),\\[2mm]
\overline{\lambda}\cdot2r^{n+1}\cdot\frac{r^{n+1}-r^{n}}{\tau }
 +(1-\overline{\lambda})\frac{\mathcal{E}_{_\mathcal{N}}[\phi^{n+1}]-\mathcal{E}_{_\mathcal{N}}[\phi^{n}]}{\tau}
 \le\left(\frac{r^{n+1}}{\sqrt{\mathcal{E}_{_\mathcal{N}}[\phi^{n}]+C}}H(\phi^{n}),\frac{\phi^{n+1}-\phi^{n}}{\tau }\right).
\end{array}\right.\label{eqremb}
\end{align}
Compared with the original weighted SAV-BE scheme \eqref{eqc3},
we replace the last equation in \eqref{eqc3} with an inequality constraint,
which not only maintains unconditional energy stability but also broadens the solution space for $\lambda$.
By applying the same decoupling strategy (see \eqref{eqc4}-\eqref{eqc6} for details),
the above constrained optimization problem can be simplified to
\begin{align}
\lambda=&\min\limits_{\overline{\lambda}\in[0,1],\,r^{n+1}}\overline{\lambda},\label{eqremc}\\[1mm]
\text{s.t.}~&~~~~
\overline{\lambda}\cdot2r^{n+1}(r^{n+1}-r^{n})
 +(1-\overline{\lambda})(\mathcal{E}_{_\mathcal{N}}[p^{n+1}+ r^{n+1}q^{n+1}]-\mathcal{E}_{_\mathcal{N}}[\phi^{n}])\nonumber\\
&\le\frac{r^{n+1}}{\sqrt{\mathcal{E}_{_\mathcal{N}}[\phi^{n}]+C}}\left(H(\phi^{n}),p^{n+1}+ r^{n+1}q^{n+1}-\phi^{n}\right).\label{eqremd}
\end{align}
\end{remark}

\subsection{Second-order weighted SAV scheme}

Based on the Crank--Nicolson (CN) formula,
we can also construct the second-order weighted SAV-CN scheme as follows:
\begin{equation}
\left\{\begin{array}{l}
\frac{\phi^{n+1}-\phi^{n}}{\tau }=-\mathcal{G}\mu^{n+\frac{1}{2}},\\[1mm]
\mu^{n+\frac{1}{2}}=\mathcal{L}\frac{\phi^{n+1}+\phi^{n}}{2}
  +\frac{r^{n+1}+r^{n}}{2\sqrt{\mathcal{E}_{_\mathcal{N}}[\phi_{*}^{n+\frac{1}{2}}]+C}}H(\phi_{*}^{n+\frac{1}{2}}),\\[2mm]
\lambda\cdot(r^{n+1}+r^{n})\cdot\frac{r^{n+1}-r^{n}}{\tau }
 +(1-\lambda)\frac{\mathcal{E}_{_\mathcal{N}}[\phi^{n+1}]-\mathcal{E}_{_\mathcal{N}}[\phi^{n}]}{\tau}
=\frac{r^{n+1}+r^{n}}{2\sqrt{\mathcal{E}_{_\mathcal{N}}[\phi_{*}^{n+\frac{1}{2}}]+C}}\left(H(\phi_{*}^{n+\frac{1}{2}}),
  \frac{\phi^{n+1}-\phi^{n}}{\tau }\right),
\end{array}\right.\label{eqd1}
\end{equation}
where $\tau >0$ is the time step size
and $\phi_{*}^{n+\frac{1}{2}}$ can be any approximation of $\phi(t^{n+\frac{1}{2}})$ with an error of $\mathcal{O}(\tau^{2})$.
For instance,
we can adopt the following backward Euler scheme with a half-time step size to obtain $\phi_{*}^{n+\frac{1}{2}}$,
\begin{equation}
\frac{\phi_{*}^{n+\frac{1}{2}}-\phi^{n}}{\tau /2}=-\mathcal{G}\left(\mathcal{L}\phi_{*}^{n+\frac{1}{2}}+H(\phi^{n})\right).\label{eqd1a}
\end{equation}

To decouple the equations for different variables,
we first eliminate $\mu^{n+\frac{1}{2}}$ by using the first two equations in \eqref{eqd1}, obtaining
\begin{equation}
\left(\mathcal{I}+\frac{\tau}{2}\mathcal{G}\mathcal{L}\right)\phi^{n+1}
=\left(\mathcal{I}-\frac{\tau}{2}\mathcal{G}\mathcal{L}\right)\phi^{n}
 -\frac{\tau\mathcal{G}H(\phi_{*}^{n+\frac{1}{2}})}{2\sqrt{\mathcal{E}_{_\mathcal{N}}[\phi_{*}^{n+\frac{1}{2}}]+C}}(r^{n+1}+r^{n}),
\label{eqd2}
\end{equation}
where $\mathcal{I}$ represents the identity operator.
Next, by decomposing $\phi^{n+1}$ from \eqref{eqd2} as
\begin{equation}
\phi^{n+1}=p_{*}^{n+1}+ r^{n+1}q_{*}^{n+1},\label{eqd3}
\end{equation}
we can derive the equations for $p_{*}^{n+1}$ and $q_{*}^{n+1}$ respectively:
\begin{equation}
\left(\mathcal{I}+\frac{\tau}{2}\mathcal{G}\mathcal{L}\right)p_{*}^{n+1}
 =\left(\mathcal{I}-\frac{\tau}{2}\mathcal{G}\mathcal{L}\right)\phi^{n}
  -\frac{\tau r^{n}\mathcal{G}H(\phi_{*}^{n+\frac{1}{2}})}{2\sqrt{\mathcal{E}_{_\mathcal{N}}[\phi_{*}^{n+\frac{1}{2}}]+C}}~~\text{and}~~
\left(\mathcal{I}+\frac{\tau}{2}\mathcal{G}\mathcal{L}\right)q_{*}^{n+1}
 =-\frac{\tau \mathcal{G}H(\phi_{*}^{n+\frac{1}{2}})}{2\sqrt{\mathcal{E}_{_\mathcal{N}}[\phi_{*}^{n+\frac{1}{2}}]+C}},\label{eqd4}
\end{equation}
in which $(\mathcal{I}+\tau\mathcal{G}\mathcal{L}/2)$ is invertible for any $\tau>0$.
Finally, substituting the solutions $p_{*}^{n+1}$ and $q_{*}^{n+1}$ of the above two equations into the third equation in \eqref{eqd1} yields
\begin{align}
0=&\,\lambda (r^{n+1}+r^{n})(r^{n+1}-r^{n})
+(1-\lambda)\left(\mathcal{E}_{_\mathcal{N}}[p_{*}^{n+1}+ r^{n+1}q_{*}^{n+1}]-\mathcal{E}_{_\mathcal{N}}[\phi^{n}]\right)\nonumber\\
&\,-\frac{r^{n+1}+r^{n}}{2\sqrt{\mathcal{E}_{_\mathcal{N}}[\phi_{*}^{n+\frac{1}{2}}]+C}}\left(H(\phi_{*}^{n+\frac{1}{2}}),
   p_{*}^{n+1}+ r^{n+1}q_{*}^{n+1}-\phi^{n}\right),\label{eqd5}
\end{align}
a nonlinear scalar equation in terms of $r^{n+1}$ which can also be solved by Newton's method with an initial guess value of $\sqrt{\mathcal{E}_{_\mathcal{N}}[\phi_{*}^{n+\frac{1}{2}}]+C}$.

Furthermore, similar results on the existence and energy stability can also be established for the second-order weighted SAV-CN scheme \eqref{eqd1}. We will omit the detailed proof, as it follows a structure akin to the first-order case.

\begin{thm}\label{thmf4}
For any given time step size $\tau>0$, the second-order weighted SAV-CN scheme \eqref{eqd1} ensures the existence of solutions when the weight coefficients are sufficiently large. More specifically, there exists a constant $\underline{\lambda}\in(0,1)$ such that the scheme always has solutions for all $\lambda\in(\underline{\lambda},1]$
\end{thm}

\begin{thm}\label{thmf5}
The second-order weighted SAV-CN scheme \eqref{eqd1} is unconditionally energy stable in the sense that
\begin{equation}
\frac{\overline{\mathcal{E}}[\phi^{n+1},r^{n+1}]-\overline{\mathcal{E}}[\phi^{n},r^{n}]}{\tau}= -(\mu^{n+1},\mathcal{G}\mu^{n+1})\le 0,\label{thd1}
\end{equation}
where
\begin{equation}
\overline{\mathcal{E}}[\phi^{n},r^{n}]=\frac{1}{2}(\phi^{n},\mathcal{L}\phi^{n})
  +\lambda((r^{n})^{2}-C)+(1-\lambda)\mathcal{E}_{_\mathcal{N}}[\phi^{n}].\label{thd2}
\end{equation}
Furthermore, if $(r^{n+1})^{2}-(r^{n})^{2}\ge \mathcal{E}_{_\mathcal{N}}[\phi^{n+1}]-\mathcal{E}_{_\mathcal{N}}[\phi^{n}]$,
the proposed scheme also adheres to the following original energy dissipation law:
\begin{equation}
\frac{\mathcal{E}[\phi^{n+1}]-\mathcal{E}[\phi^{n}]}{\tau}\le -(\mu^{n+1},\mathcal{G}\mu^{n+1})\le 0,
\quad\text{where}\quad
\mathcal{E}[\phi^{n}]=\frac{1}{2}(\phi^{n},\mathcal{L}\phi^{n})+\mathcal{E}_{_\mathcal{N}}[\phi^{n}].\label{thd3}
\end{equation}
\end{thm}

\begin{thm}\label{thmf6}
If the weight coefficient $\lambda=\lambda_{\star}\in[0,1]$ yields a solution $(\phi_{\star}^{n+1},\mu_{\star}^{n+1},r_{\star}^{n+1})$ for the weighted SAV-CN scheme \eqref{eqd1}, then

(i) when $(r^{n+1})^{2}-(r^{n})^{2}\ge \mathcal{E}_{_\mathcal{N}}[\phi_{\star}^{n+1}]-\mathcal{E}_{_\mathcal{N}}[\phi^{n}]$,
the scheme \eqref{eqd1} always has solutions for any $\lambda\in[0,\lambda_{\star}]$;

(ii) when $(r^{n+1})^{2}-(r^{n})^{2}<\mathcal{E}_{_\mathcal{N}}[\phi_{\star}^{n+1}]-\mathcal{E}_{_\mathcal{N}}[\phi^{n}]$,
the scheme \eqref{eqd1} always has solutions for any $\lambda\in[\lambda_{\star},1]$.
\end{thm}

\begin{algorithm}
\caption{The second-order weighted SAV-CN scheme \eqref{eqd1}.}\label{alg2}
Update the unknown variables $(\phi^{n+1},r^{n+1})$ by the known ones $(\phi^{n},r^{n})$.
\begin{enumerate}
  \item Solve $p_{*}^{n+1}$ and $q_{*}^{n+1}$ form \eqref{eqd4}.
  \item Set $\lambda =0$ and solve Eq.~\eqref{eqd5}.\\
  \textbf{if} Eq.~\eqref{eqc7} is solvable \textbf{then} \\
  $~~~~$ Output the final values of $r^{n+1}$ and $\lambda=0$.\\
  \textbf{else}  \\
  $~~~~$ Set $a=0, b=1$ and provide the tolerance $tol$.\\
  $~~~~$ \textbf{while} $b-a\ge tol$ \textbf{do}\\
  $~~~~~~~~~$ Set $\lambda=\frac{b+a}{2}$ and solve Eq.~\eqref{eqd5}.\\
  $~~~~~~~~~$ \textbf{if} Eq.~\eqref{eqc7} is solvable \textbf{then}\\
  $~~~~~~~~~~~~~$ Set $a=a$ and $b=\frac{b+a}{2}$.\\
  $~~~~~~~~~$ \textbf{else}\\
  $~~~~~~~~~~~~~$ Set $a=\frac{b+a}{2}$ and $b=b$.\\
  $~~~~~~~~~$ \textbf{end if}\\
  $~~~~$ \textbf{end while}\\
  $~~~~$ Output the final values of $r^{n+1}$ and $\lambda$.\\
  \textbf{end if}
  \item Update $\phi^{n+1}$ by \eqref{eqd3}.
\end{enumerate}
\end{algorithm}

Based on the above theoretical results,
we develop Algorithm \ref{alg2} for the second-order SAV-CN scheme \eqref{eqd1}.
It is worth noting that, for the gradient flow equation \eqref{eqc1}, alternative time-discretization methods such as linear multi-step and Runge--Kutta formulas can also be employed to develop unconditionally energy stable schemes \cite{Huang23,Akrivis19,GongY20C,Feng21}.
Furthermore, since the proofs of unconditional energy stability for weighted SAV schemes are variational in nature,
these results can be readily extended to the fully discrete setting. Various methods, including the Galerkin finite element method \cite{ChenH20}, finite difference method \cite{Huang23b}, and spectral method \cite{Huang23}, can be employed for spatial discretization.
On the other hand, since the key to the well-posedness lies in solving the nonlinear scalar equations \eqref{eqc7} and \eqref{eqd5},
it follows that the spatial discretization does not impact the existence of their solutions.

\section{Numerical results}
\label{sec4}

In this section, various numerical examples will be presented to verify the effectiveness and energy stability of our weighted SAV method.
Employing periodic boundary conditions,
we shall utilize the Fourier spectral method for spatial discretization and the composite trapezoidal rule for numerical integration.
The tolerance for the weighting coefficient $\lambda$ is set to $tol=10^{-8}$.
The following widely used Ginzburg--Landau type free energy~\cite{ShenJ18,Huang19b,JiangM22} will be considered:
\begin{equation}
\mathcal{E}[\phi]=\int_{\Omega}\left[\frac{\varepsilon^{2}}{2}|\nabla\phi|^{2}+\frac{1}{4}(1-\phi^{2})^{2}\right]\mathrm{d}\bm{x}.\label{eqe2}
\end{equation}
Here $0<\varepsilon\ll 1$ is a small parameter to measure the thickness of the interface.
As is customary~\cite{ShenJ18,JiangM22,Huang24}, neglecting irrelevant constant terms, the above free energy is divided into two parts,
$\mathcal{L}=-\varepsilon^{2}\Delta+\gamma$ and $\mathcal{E}_{_\mathcal{N}}[\phi]=\int_{\Omega}F(\phi)\mathrm{d}\bm{x}$,
where $F(\phi)=\frac{1}{4}(\phi^{2}-1-\gamma)^{2}$ with a constant $\gamma\ge 0$.
Furthermore, to ensure that $\mathcal{E}_{_\mathcal{N}}[\phi]$ remains theoretically bounded, we apply the following truncation:
\begin{equation}
F_{\delta}(\phi)=\left\{\begin{array}{ll}
(a\phi+b)\mathrm{e}^{-\phi}+c,&\phi>\delta,\\[1mm]
F(\phi),& -\delta\le\phi\le\delta,\\[1mm]
(-a\phi+b)\mathrm{e}^{\phi}+c,&\phi<-\delta,
\end{array}\right.\label{eqe3}
\end{equation}
where
\begin{equation}
\left\{\begin{array}{l}
a=-\left[\delta^{3}+3\delta^{2}-(1+\gamma)\delta-(1+\gamma)\right]\mathrm{e}^{\delta},\\[1mm]
b=\left[\delta^{4}+\delta^{3}-(4+\gamma)\delta^{2}+(1+\gamma)\delta+(1+\gamma)\right]\mathrm{e}^{\delta},\\[1mm]
c=\frac{1}{4}\delta^{4}+2\delta^{3}+\frac{5-\gamma}{2}\delta^{2}-2(1+\gamma)\delta+\frac{(1+\gamma)(\gamma-3)}{4},
\end{array}\right.\label{eqe4}
\end{equation}
with $\delta$ is a parameter.
It is easy to verify that the above truncation ensures $F_{\delta}(\phi)$ remains an even function and $F_{\delta}(\phi)\in C^{2}(\mathbb{R})$.
In numerical simulations, we set $0\le\gamma\le 4$ and $\delta=5$ in order to ensure that
the stationary points of $F_{\delta}(\phi)$ precisely coincide with those of $F(\phi)$, while preserving its monotonicity.

Applying the $H^{-\nu}$-gradient flow (where $\mathcal{G}=(-\Delta)^{\nu},\,0<\nu<1$) to the free energy~\eqref{eqe2},
yields the space-fractional Cahn--Hilliard equation~\cite{Ainsworth17,TangT19}.
Specifically, when $\nu=0$ (i.e. $L^{2}$-gradient flow) and $\nu=1$ (i.e., $H^{-1}$-gradient flow),
they correspond to the commonly-used Allen--Cahn equation and Cahn--Hilliard equation, respectively.
Furthermore, when $0<\nu\le 1$, the space-fractional Cahn--Hilliard equation also satisfies the following conservation law of mass,
\begin{equation}
\frac{\mathrm{d}m(t)}{\mathrm{d}t}=0,\quad\text{where}\quad m(t)=\int_{\Omega}\phi(\bm{x},t)\mathrm{d}\bm{x}.\label{eqe5}
\end{equation}
It is straightforward to prove that our proposed weighted SAV scheme preserves this property discretely.

\subsection{Convergence rate}

To begin with, we examine the convergence rate of the proposed weighted SAV scheme
by using the Cahn--Hilliard equation with $\varepsilon=0.1$ and $\gamma=4$.
For this purpose, we consider the initial conditions~\cite{ShenJ18,LiuZ20,Huang23} in the following form:
\begin{equation}
\phi(x,y,0)=0.05\sin(x)\cos(y).\label{eqe6}
\end{equation}
Moreover, the computational region $\Omega=[0,2\pi)^{2}$ is discretized into $128\times 128$ Fourier modes.

\begin{table}[!htp]
\renewcommand\arraystretch{1.4}
\centering
\caption{Temporal errors in the $l^{2}$-norm and $l^{\infty}$-norm at $T=0.5$ for the weighted SAV-BE scheme \eqref{eqc3}.}\label{tabrate1}
\tabcolsep 0.055in
\begin{tabular}[c]
{ccccccccccccc} \hline
&&\multicolumn{5}{c}{$\lambda=0.5$}&&\multicolumn{5}{c}{$\lambda=\lambda_{\min}$}\\
\cline{3-7} \cline{9-13}
$\tau $ &&$\|\phi_{\tau}-\phi_{\tau/2}\|_{2}$ &Rate &&$\|\phi_{\tau}-\phi_{\tau/2}\|_{\infty}$ &Rate &&$\|\phi_{\tau}-\phi_{\tau/2}\|_{2}$ &Rate &&$\|\phi_{\tau}-\phi_{\tau/2}\|_{\infty}$ &Rate\\ \hline
$0.5^{2}\times10^{-3}$ &&7.0243e-3&        &&3.1134e-3&        &&7.0240e-3&        &&3.1132e-3&        \\
$0.5^{3}\times10^{-3}$ &&4.0641e-3& 0.7894 &&1.8213e-3& 0.7735 &&4.0639e-3& 0.7894 &&1.8212e-3& 0.7735 \\
$0.5^{4}\times10^{-3}$ &&2.2012e-3& 0.8846 &&9.9072e-4& 0.8784 &&2.2011e-3& 0.8846 &&9.9066e-4& 0.8784 \\
$0.5^{5}\times10^{-3}$ &&1.1478e-3& 0.9395 &&5.1750e-4& 0.9369 &&1.1477e-3& 0.9395 &&5.1746e-4& 0.9369 \\
$0.5^{6}\times10^{-3}$ &&5.8634e-4& 0.9690 &&2.6458e-4& 0.9679 &&5.8631e-4& 0.9690 &&2.6456e-4& 0.9679 \\
$0.5^{7}\times10^{-3}$ &&2.9637e-4& 0.9843 &&1.3378e-4& 0.9838 &&2.9636e-4& 0.9843 &&1.3377e-4& 0.9838 \\
\hline
\end{tabular}
\end{table}

\begin{table}[!htp]
\renewcommand\arraystretch{1.4}
\centering
\caption{Temporal errors in the $l^{2}$-norm and $l^{\infty}$-norm at $T=0.5$ for the weighted SAV-CN scheme \eqref{eqd1}.}\label{tabrate2}
\tabcolsep 0.055in
\begin{tabular}[c]
{ccccccccccccc} \hline
&&\multicolumn{5}{c}{$\lambda=0.5$}&&\multicolumn{5}{c}{$\lambda=\lambda_{\min}$}\\
\cline{3-7} \cline{9-13}
$\tau $ &&$\|\phi_{\tau}-\phi_{\tau/2}\|_{2}$ &Rate &&$\|\phi_{\tau}-\phi_{\tau/2}\|_{\infty}$ &Rate &&$\|\phi_{\tau}-\phi_{\tau/2}\|_{2}$ &Rate &&$\|\phi_{\tau}-\phi_{\tau/2}\|_{\infty}$ &Rate\\ \hline
$0.5^{2}\times10^{-3}$ &&1.5943e-4&        &&7.0923e-5&        &&1.5906e-4&        &&7.0719e-5&        \\
$0.5^{3}\times10^{-3}$ &&4.1174e-5& 1.9531 &&1.8313e-5& 1.9534 &&4.0986e-5& 1.9563 &&1.8209e-5& 1.9575 \\
$0.5^{4}\times10^{-3}$ &&1.0455e-5& 1.9775 &&4.6498e-6& 1.9776 &&1.0361e-5& 1.9840 &&4.5972e-6& 1.9858 \\
$0.5^{5}\times10^{-3}$ &&2.6338e-6& 1.9890 &&1.1711e-6& 1.9893 &&2.5865e-6& 2.0021 &&1.1447e-6& 2.0058 \\
$0.5^{6}\times10^{-3}$ &&6.6070e-7& 1.9951 &&2.9351e-7& 1.9965 &&6.3717e-7& 2.0212 &&2.8027e-7& 2.0301 \\
$0.5^{7}\times10^{-3}$ &&1.6511e-7& 2.0005 &&7.2839e-8& 2.0106 &&1.5349e-7& 2.0536 &&6.6150e-8& 2.0830 \\
\hline
\end{tabular}
\end{table}

Tables \ref{tabrate1} and \ref{tabrate2} respectively present the numerical errors and convergence rates of
the weighted SAV-BE scheme \eqref{eqc3} and the weighted SAV-CN scheme \eqref{eqd1} with respect to different time step sizes.
For both the nonlinear energy-based type (corresponding to $\lambda=1$)
and the Lagrange multiplier type (corresponding to $\lambda=0$) SAV schemes,
their convergence rates have been extensively verified~\cite{ShenJ18,ChengQ20,Huang23}.
Here, we further consider a particular intermediate state ($\lambda=0.5$) and the state using the minimum weighting coefficients ($\lambda=\lambda_{\min}$), which illustrate that our scheme still exhibits the first and second order of convergence rates.

\subsection{Dynamic evolution in 2D}

Next, we consider initial conditions with the following distribution
\begin{equation}
\phi(\bm{x},t)|_{t=0}=\tanh\left(\frac{\mathrm{dist}(\bm{x},\Gamma)}{\sqrt{2}\varepsilon}\right),\label{eqini}
\end{equation}
where $\mathrm{dist}(\bm{x},\Gamma)$ represents the signed distance from point $\bm{x}$ to some curve $\Gamma$. For example, $\Gamma$ can be selected as the boundary of a cross-shaped region defined by
\begin{equation}
\Omega_{0}=\{(x,y)|(x,y)\in [-0.75,0.75]\times[-0.25,0.25]\}\cup \{(x,y)|(x,y)\in [-0.25,0.25]\times[-0.75,0.75\}.\label{eqcurve1a}
\end{equation}
In this subsection, we explore the dynamics governed by the Cahn--Hilliard equation with $\varepsilon=0.01$
for the following closed curves~\cite{JiangW24,Pei23},
\begin{align}
&\text{Curve I: \quad boundary of the cross-shaped region $\Omega_0$},\label{eqcurve1}\\
&\text{Curve II:}\quad\left\{\begin{array}{l}
x=\cos\theta,\\[1mm]
y=2\sin\theta-1.9\sin^{3}\theta,
\end{array}\right.\label{eqcurve2}\\
&\text{Curve III:}\quad\left\{\begin{array}{l}
x=0.25(3\cos\theta+\cos(3\theta)),\\[1mm]
y=0.25(3\sin\theta-\sin(3\theta)),
\end{array}\right.\label{eqcurve3}\\
&\text{Curve IV:}\quad\left\{\begin{array}{l}
x=0.5\cos\theta,\\[1mm]
y=0.25\sin\theta+0.5\sin(\cos\theta)+0.5(0.2+\sin\theta\sin^{2}(3\theta))\sin\theta,
\end{array}\right.\label{eqcurve4}
\end{align}
where $\theta\in[0,2\pi)$.
The computational domain $\Omega$ is discretized with $256\times 256$ modes.

\begin{figure}[htp]
\centering
\includegraphics[width=15cm]{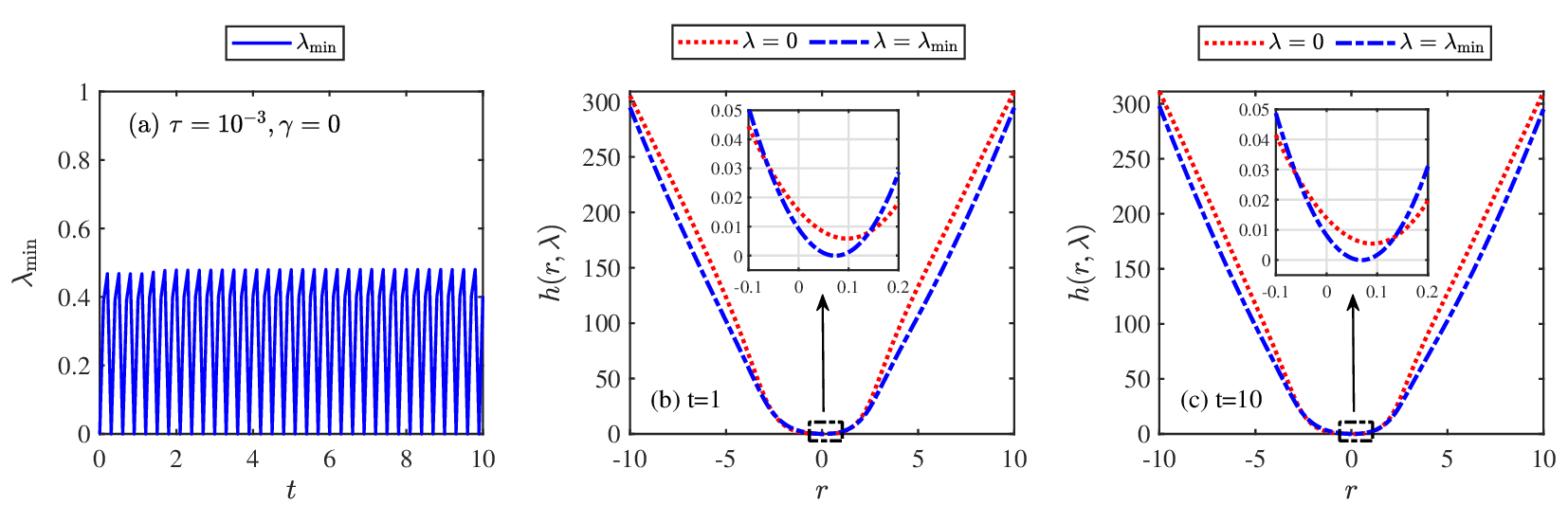}
\caption{(a): The evolution of weighting coefficient $\lambda_{\min}$ of curve I simulated using the weighted SAV-BE scheme \eqref{eqc3}
         with $\tau=10^{-3}$ and $\gamma=0$.
         (b)-(c): The curve of the left-hand side of nonlinear scalar equation \eqref{eqc7} (or $h(r,\lambda)$ as defined in \eqref{thmh1})
         with respect to $r$ at $t=1$ and $t=10$.}
\label{figc1a}
\end{figure}

\begin{figure}[htp]
\centering
\includegraphics[width=15cm]{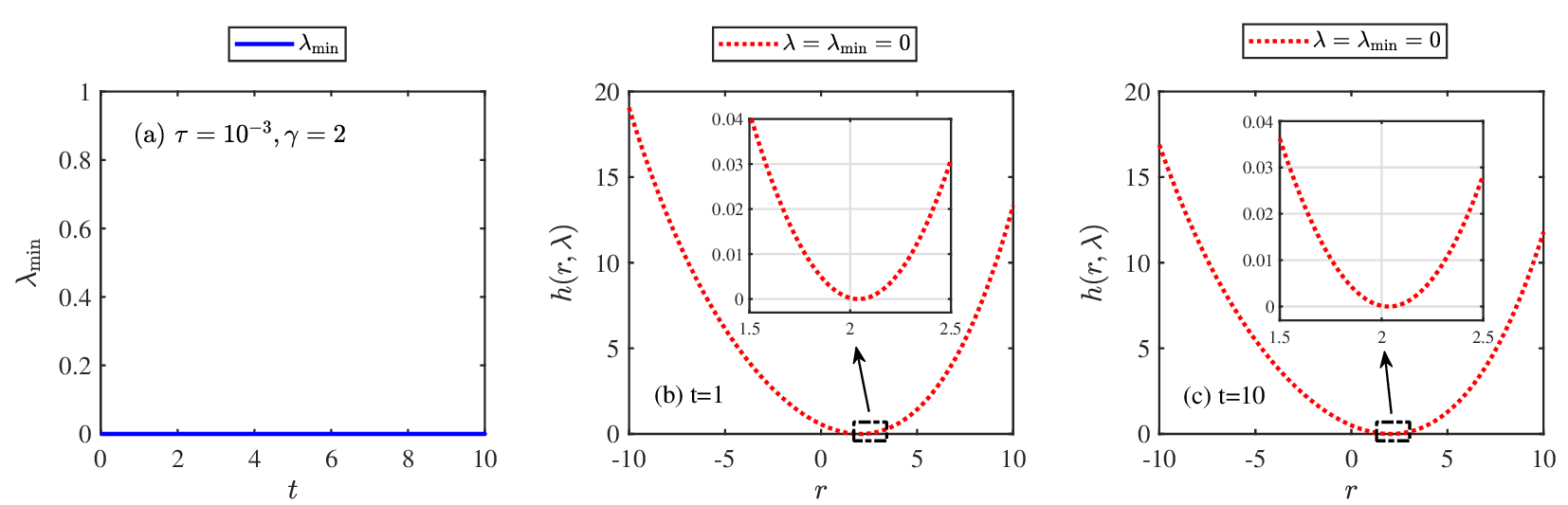}
\caption{(a): The evolution of weighting coefficient $\lambda_{\min}$ of curve I simulated using the weighted SAV-BE scheme \eqref{eqc3}
         with $\tau=10^{-3}$ and $\gamma=2$.
         (b)-(c): The curve of the left-hand side of nonlinear scalar equation \eqref{eqc7} (or $h(r,\lambda)$ as defined in \eqref{thmh1})
         with respect to $r$ at $t=1$ and $t=10$.}
\label{figc1b}
\end{figure}

We first examine the behavior of the nonlinear scalar equation \eqref{eqc7}. Taking the evolution of curve I
with parameters $\tau=10^{-3}$ and $\gamma=0$ as an example, we first plot the weight coefficient $\lambda_{\min}$ as a function of time in Fig.~\ref{figc1a}(a), which demonstrates a frequent oscillation between $0$ and $0.5$ during the evolution, and it implies that the Lagrange multiplier-based approach (i.e., $\lambda\equiv 0$) has encountered severe troubles with no solution. Furthermore, at times $t=1, 10$, we depict the function curve $h(r,\lambda)$ in Figs.~\ref{figc1a}(b) and \ref{figc1a}(c) respectively, which clearly illustrate that the nonlinear scalar equation \eqref{eqc7} has no solution under the Lagrange multiplier framework ($\lambda=0$), while with $\lambda_{\min}$ obtained by our weighted approach $h(r,\lambda_{\min})=0$ does own a solution. On the other hand, it is widely recognized that the parameter $\gamma$ typically serves as a stabilizer in numerical schemes; larger values of $\gamma$ enhance stability but can also significantly increase numerical errors.
Interestingly, a comparison between the results in Figs.~\ref{figc1a} and \ref{figc1b} shows that the stabilizer $\gamma$
also affects the existence of solutions to the nonlinear scalar equation \eqref{eqc7}.

\begin{figure}[htbp!]
\centering
\includegraphics[width=15cm]{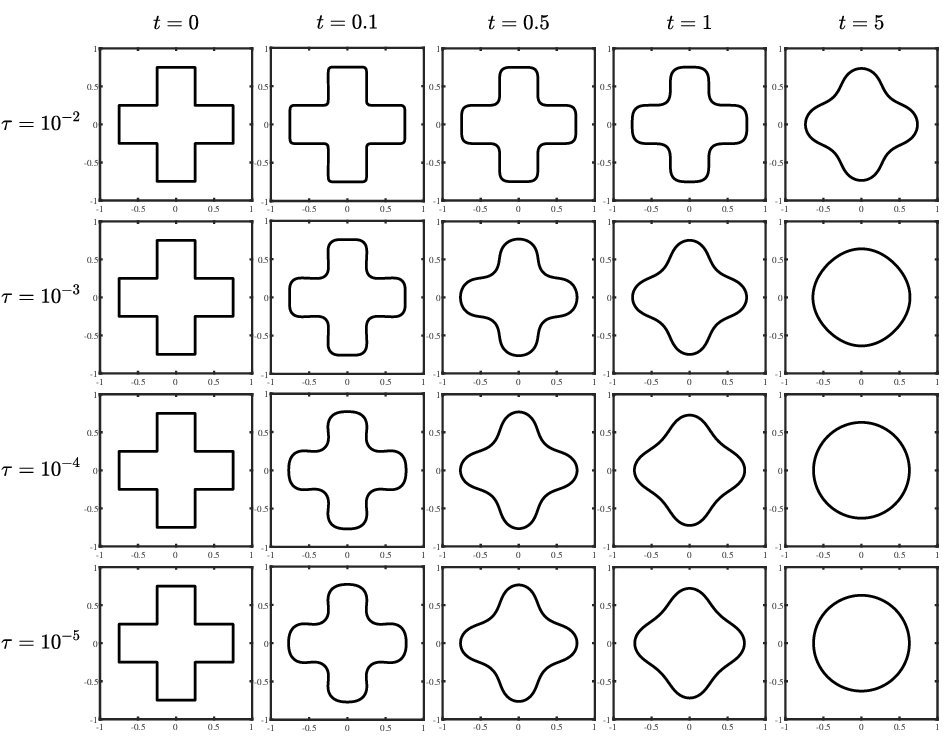}
\caption{Several snapshots of curve I simulated using the weighted SAV-BE scheme \eqref{eqc3} with different $\tau$, respectively, where $\gamma=2$.}
\label{figc1c}
\end{figure}

\begin{figure}[htbp!]
\centering
\includegraphics[width=15cm]{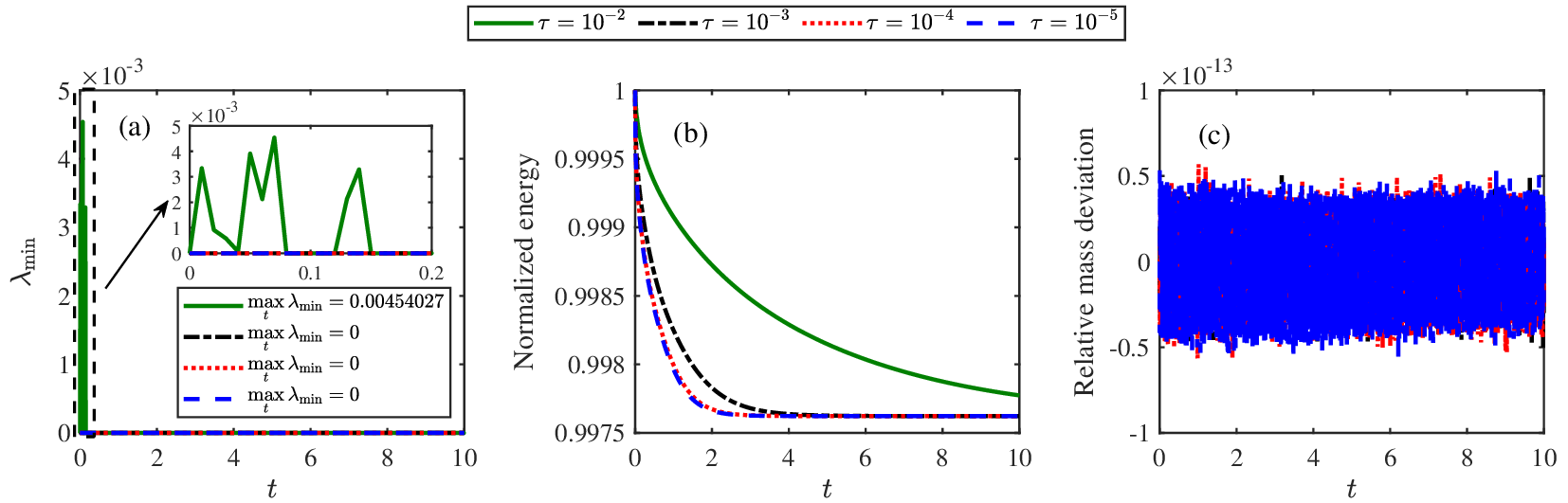}
\caption{The evolution of (a) weighting coefficient $\lambda_{\min}$, (b) normalized energy, and (c) relative mass deviation
          of curve I simulated using the weighted SAV-BE scheme \eqref{eqc3} with different time step sizes, respectively, where $\gamma=2$.}
\label{figc1d}
\end{figure}

Fig.~\ref{figc1c} plots the evolution of curve I under different time step sizes.
It is evident that with larger time step sizes, the dynamical evolution exhibits noticeable lag compared to smaller time step sizes,
consistent with the numerical findings in~\cite{Huang23}.
Moreover, the results in Fig.~\ref{figc1d}(a) clearly show that with larger time step sizes,
the weighting coefficients may be non-zero,
thereby numerically confirming the theoretical findings in~\cite{Onuma24}.
This result is corroborated by the energy evolution depicted in Fig.~\ref{figc1d}(b).
In addition, the conservation of mass in the Cahn--Hilliard equation is validated at the discrete level, as shown in Fig.~\ref{figc1d}(c).

\begin{figure}[htbp!]
\centering
\includegraphics[width=15cm]{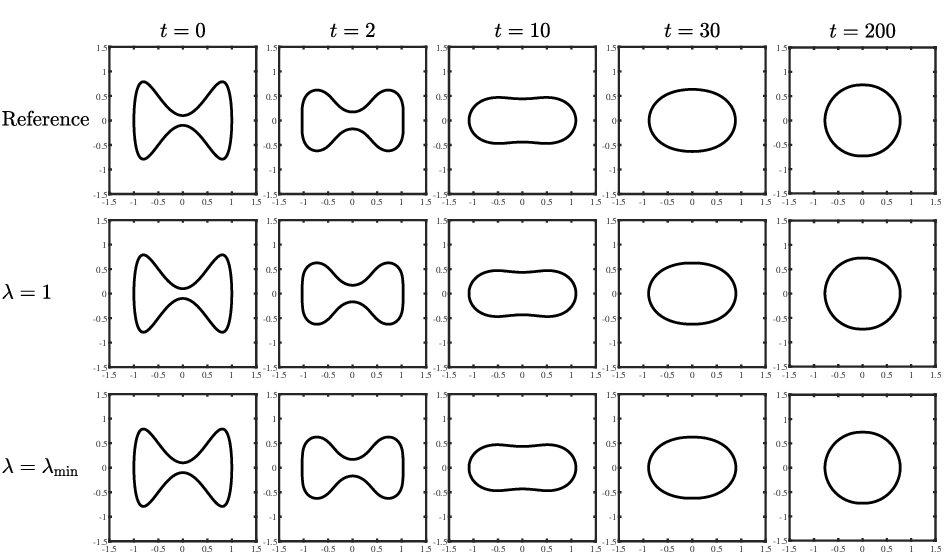}
\caption{Several snapshots of curve II simulated using the SAV-CN scheme (second row)
         and the weighted SAV-CN scheme \eqref{eqd1} (third row), respectively, where $\tau=10^{-3}$ and $\gamma=2$.}
\label{figc2a}
\end{figure}

\begin{figure}[htbp!]
\centering
\includegraphics[width=15cm]{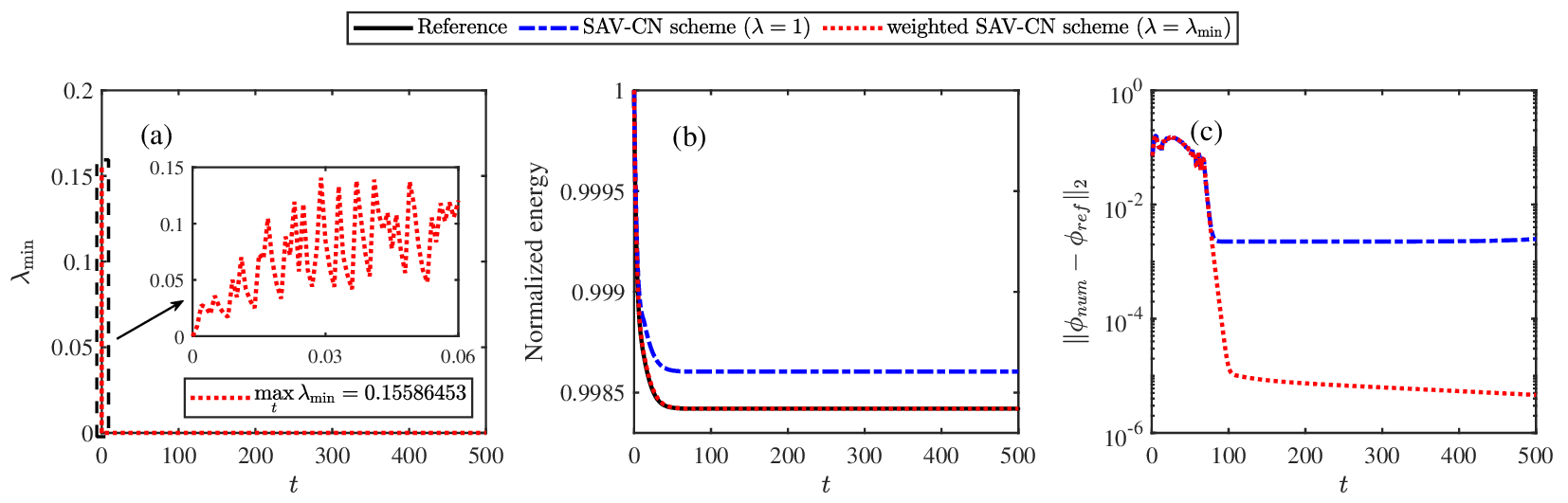}
\caption{The evolution of (a) weighting coefficient $\lambda_{\min}$, (b) normalized energy, and (c) error $\|\phi_{num}-\phi_{ref}\|_{2}$
         of curve II simulated using the SAV-CN scheme and the weighted SAV-CN scheme \eqref{eqd1}, respectively,
         where $\tau=10^{-3}$ and $\gamma=2$.}
\label{figc2b}
\end{figure}

Moreover, we compare the energy and solution behaviors of the nonlinear energy-based SAV-CN scheme (i.e., $\lambda\equiv 1$) and the weighted SAV-CN scheme (i.e., $\lambda=\lambda_{\min}$). Fig.~\ref{figc2a} plots the evolution process of curve II, which indicates that there are no significant visual differences in the evolution process by using different schemes. Nevertheless, by using the solution obtained from the Lagrange multiplier approach (i.e., $\lambda\equiv 0$) with $\tau=10^{-6}$ and $\gamma=2$ as the reference solution, we can observe a notable energy deviation between the nonlinear energy-based SAV-CN scheme and the reference energy in Fig.~\ref{figc2b}(b),
while the energy of weighted strategy exhibits no observable differences. These findings align with our theoretical expectations. Regarding the solution, Fig.~\ref{figc2b}(c) illustrates the $l^{2}$-error for various methods. It clearly shows that our method outperforms the nonlinear energy-based SAV-CN scheme in terms of accuracy as time progresses, achieving a significantly higher level of precision.

\begin{figure}[htp]
\centering
\includegraphics[width=15cm]{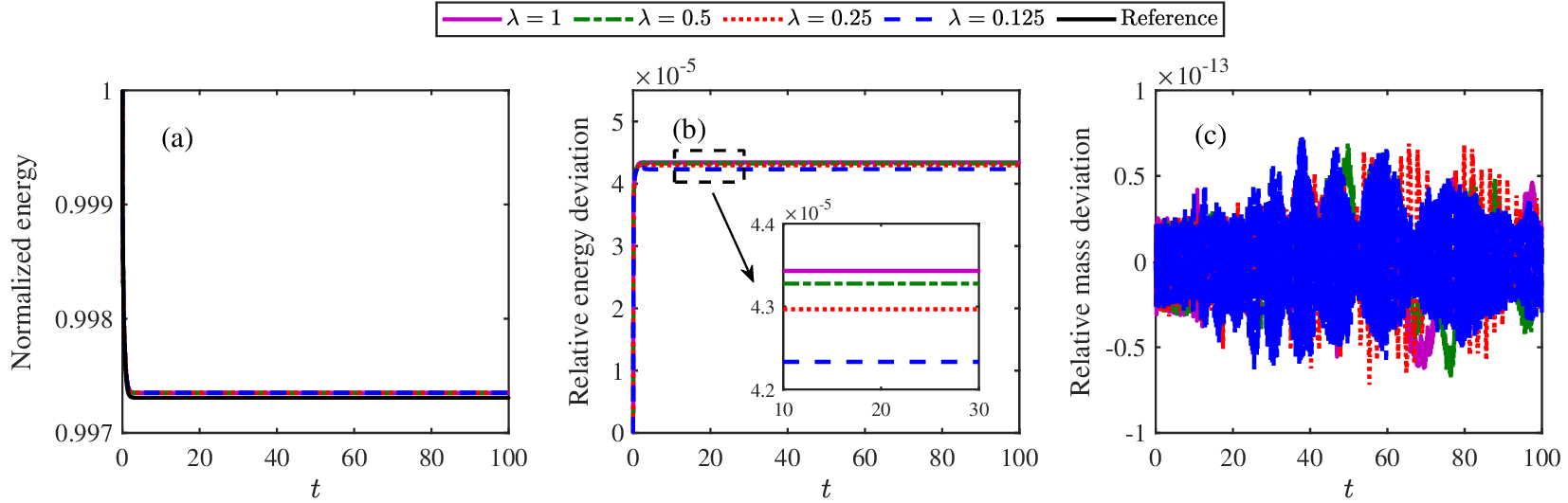}
\caption{The evolution of (a) normalized energy, (b) relative energy deviation, (c) relative mass deviation of curve III simulated
         using the weighted SAV-CN scheme \eqref{eqd1} with different weighting coefficients, respectively, where $\tau=10^{-4}$ and $\gamma=1$.}
\label{figc3b}
\end{figure}

Following this, we further numerically explore the convergence of the modified energy $\overline{\mathcal{E}}[\phi,r]$
(as defined by \eqref{thd2}) for curve III, focusing on the effects of different fixed weight coefficients $\lambda$. To ensure the solution existence of the weighted SAV-CN scheme \eqref{eqd1} under different weights, we set $\tau=10^{-4}$ and $\gamma=1$. The reference solution is obtained through the Lagrange multiplier scheme (i.e., $\lambda\equiv 0$) with $\tau=10^{-6}$ and $\gamma=1$.
As shown in Fig.~\ref{figc3b}, with the weighting coefficient $\lambda$ decreasing,
the modified energy $\overline{\mathcal{E}}[\phi^{n},r^{n}]$ gradually converges towards the reference energy $\mathcal{E}[\phi^{n}]$ at the discrete level. These findings agree with our theoretical expectations.
The mass conservation is also verified in Fig.~\ref{figc3b}(c).
For a more precise quantitative comparison, it is encouraged to refer to Table \ref{tabrate3}.
This table show that as the termination time $T$ increases, the modified energy $\overline{\mathcal{E}}[\phi,r]$
exhibits a power-law close to $-1$ with respect to the weight coefficient $\lambda$.
This novel finding warrants further exploration in our future work.

\begin{table}[!htp]
\renewcommand\arraystretch{1.4}
\centering
\caption{Energy errors $\|\overline{\mathcal{E}}_{\lambda}-\overline{\mathcal{E}}_{\lambda/2}\|_{2}
:=\sqrt{\int_{0}^{T}[\overline{\mathcal{E}}_{\lambda}(t)-\overline{\mathcal{E}}_{\lambda/2}(t)]^{2}\mathrm{d}t}$
of curve III with different termination times $T$.}\label{tabrate3}
\tabcolsep 0.065in
\begin{tabular}[c]
{cccccccccccccccc} \hline
&&\multicolumn{2}{c}{$T=0.1$}&&\multicolumn{2}{c}{$T=1$}&&\multicolumn{2}{c}{$T=10$}&&\multicolumn{2}{c}{$T=100$}\\
\cline{3-4} \cline{6-7} \cline{9-10} \cline{12-13}
$\lambda$ &&$\|\overline{\mathcal{E}}_{\lambda}-\overline{\mathcal{E}}_{\lambda/2}\|_{2}$ &Rate &&
$\|\overline{\mathcal{E}}_{\lambda}-\overline{\mathcal{E}}_{\lambda/2}\|_{2}$ &Rate &&
$\|\overline{\mathcal{E}}_{\lambda}-\overline{\mathcal{E}}_{\lambda/2}\|_{2}$ &Rate &&
$\|\overline{\mathcal{E}}_{\lambda}-\overline{\mathcal{E}}_{\lambda/2}\|_{2}$ &Rate\\ \hline
$1$       &&1.1972e-7&         &&4.1890e-7&         &&3.0500e-6&         &&1.0339e-5&        \\
$0.5$     &&2.3596e-7& -0.9789 &&8.2708e-7& -0.9814 &&6.1710e-6& -1.0167 &&2.0916e-5& -1.0165 \\
$0.5^{2}$ &&4.5782e-7& -0.9562 &&1.6152e-6& -0.9656 &&1.2601e-5& -1.0300 &&4.2684e-5& -1.0291 \\
$0.5^{3}$ &&8.5812e-7& -0.9064 &&3.1081e-6& -0.9443 &&2.6039e-5& -1.0471 &&8.8015e-5& -1.0441 \\
$0.5^{4}$ &&1.4806e-6& -0.7870 &&5.9886e-6& -0.9462 &&5.3865e-5& -1.0487 &&1.8094e-4& -1.0397 \\
$0.5^{5}$ &&2.0293e-6& -0.4548 &&1.2349e-5& -1.0441 &&1.0504e-4& -0.9636 &&3.4781e-4& -0.9428 \\
\hline
\end{tabular}
\end{table}


\begin{figure}[htp]
\centering
\includegraphics[width=15cm]{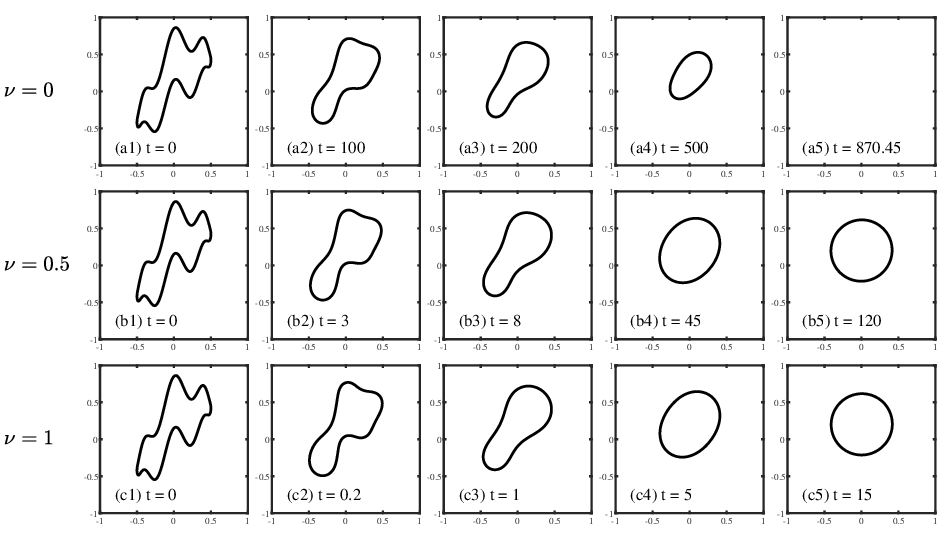}
\caption{Several snapshots of curve IV simulated using the weighted SAV-CN scheme \eqref{eqd1} with different $\nu$, respectively,
         where $\tau=10^{-3}$ and $\gamma=4$.}
\label{figc3}
\end{figure}

\begin{figure}[htp]
\centering
\includegraphics[width=15cm]{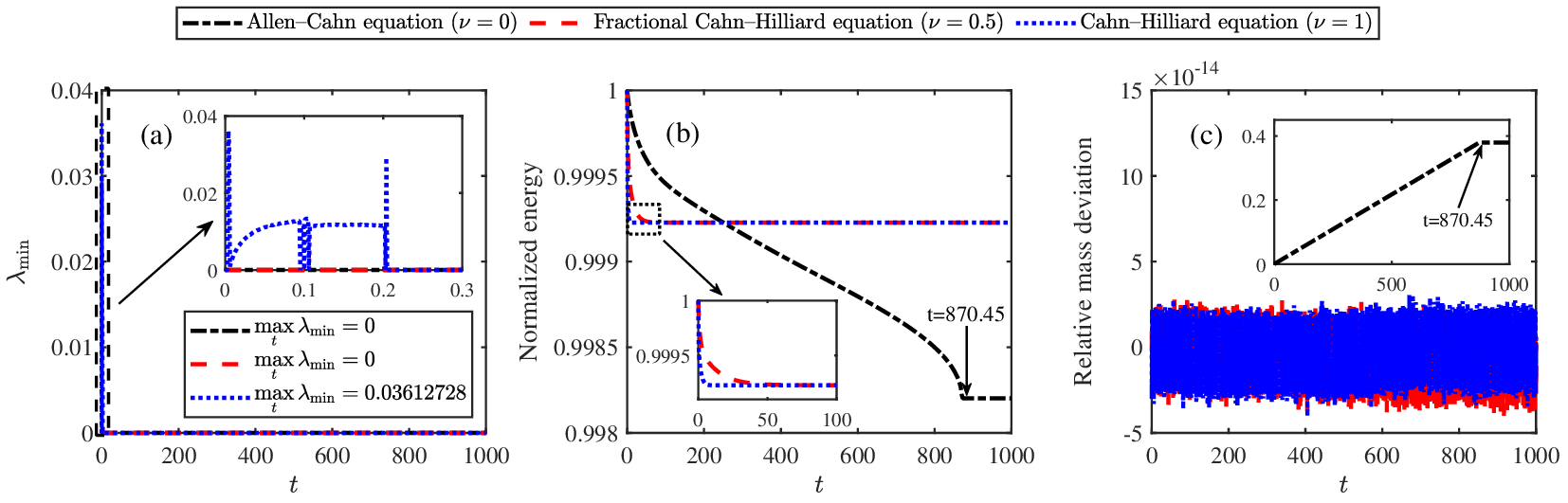}
\caption{The evolution of (a) weighting coefficient $\lambda_{\min}$, (b) normalized energy, and (c) relative mass deviation of curve IV
         simulated using the weighted SAV-CN scheme \eqref{eqd1} with different $\nu$, respectively, where $\tau=10^{-3}$ and $\gamma=4$.}
\label{figt3}
\end{figure}

Finally, we examine the evolution of curve IV under the Allen--Cahn equation,
the spatial-fractional Cahn--Hilliard equation with $\nu=0.5$, and the Cahn--Hilliard equation.
The Allen--Cahn equation is known for describing interface evolution governed by mean curvature flow~\cite{DuQ20},
where the main geometric feature of this evolution is the contraction and eventual disappearance of closed curves/surfaces.
In contrast, for the spatial-fractional Cahn--Hilliard equation ($0<\nu\le 1$), constrained by mass conservation,
its evolution is typically dominated by surface diffusion flow (often including various forms of bulk diffusion as well).
Under isotropic interfacial energy, closed curves/surfaces typically evolve towards circular/spherical shapes,
as numerically validated in Fig.~\ref{figc3}.
It is obvious that the fractional-order derivative $\nu$ significantly influences the evolution speed of the interface,
consistent with the results given in~\cite{Huang23}.
In addition, larger fractional-order derivatives $\nu$ result in higher spatial derivatives,
thereby increasing the complexity of the equation and consequently influencing the values of the weighting coefficient $\lambda$.
This phenomenon is also validated in Fig.~\ref{figt3}(a).
Due to the lack of mass conservation constraints in the Allen--Cahn equation,
the system tends to minimize its free energy globally, often resulting in trivial solutions (i.e., $\phi\equiv -1$).
In contrast, the Cahn--Hilliard equation,  which incorporates a mass conservation constraint, typically evolves towards non-trivial solutions.
This difference is also reflected in the energy curves shown in Fig.~\ref{figt3}(b),
where the curve corresponding to the Allen--Cahn equation shows lower energy levels in its stable state.
Furthermore, interfaces undergo significant energy changes during phase transitions
(corresponding to the complete disappearance of the closed curve at $t=870.45$ under the Allen--Cahn equation case).
This property is also observed in many evolutionary processes (please refer to \cite{BaoW17,Huang19a,Huang23b}).

\subsection{Dynamic evolution in 3D}

As the last example, we simulate the evolution of a torus~\cite{DuQ06,Yoon21} in three-dimensional space using the Cahn--Hilliard equation,
with a major radius $R=0.6$ and a minor radius $r=0.3$, while maintaining the initial surface profile described in \eqref{eqini}.
We set the phase-field parameters $\varepsilon=0.02$ and $\gamma=4$, the time step size $\tau=10^{-4}$,
and let the spatial domain $\Omega=[-1,1)^{3}$ be discretized into $128^{3}$ Fourier modes.

\begin{figure}[!htp]
\centering
\includegraphics[width=15cm]{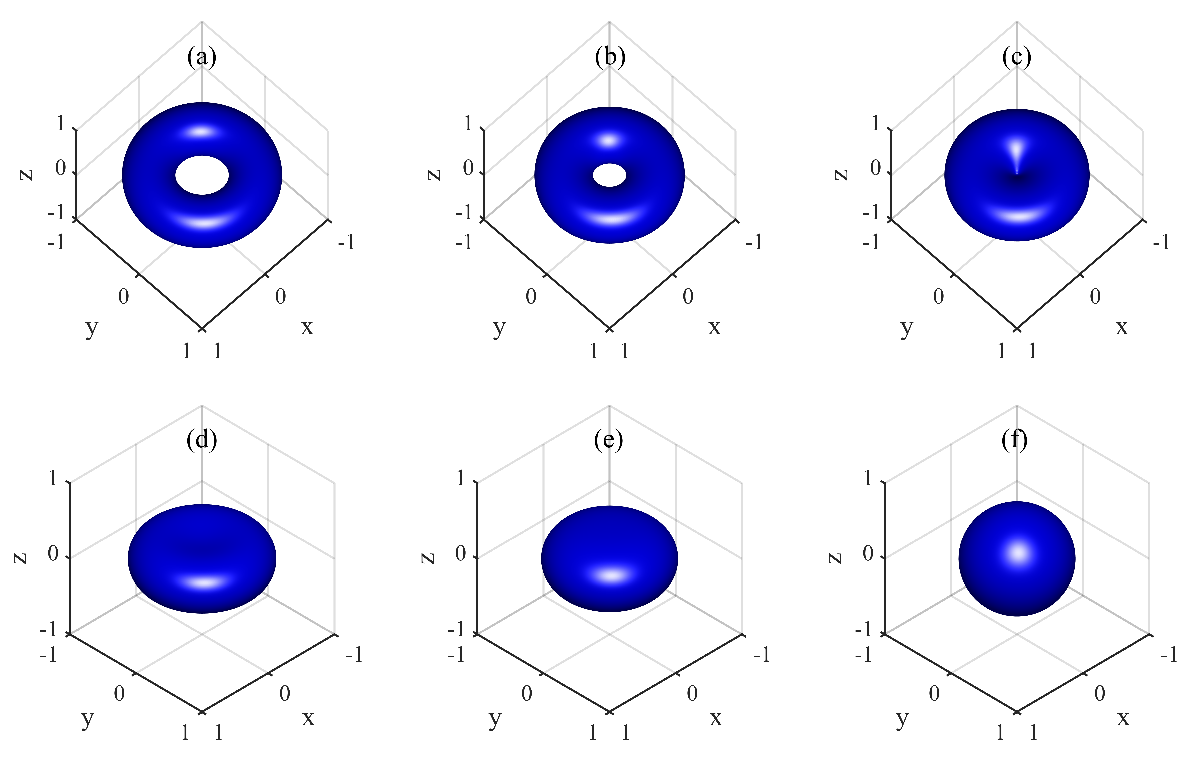}
\caption{Morphological evolution of an initial torus at six times:
         (a) $t=0$, (b) $t=1$, (c) $t=1.583$, (d) $t=2$, (e) $t=3$, (f) $t=10$ simulated using the weighted SAV-CN scheme \eqref{eqd1}.}
\label{figc4}
\end{figure}

\begin{figure}[htp]
\centering
\includegraphics[width=15cm]{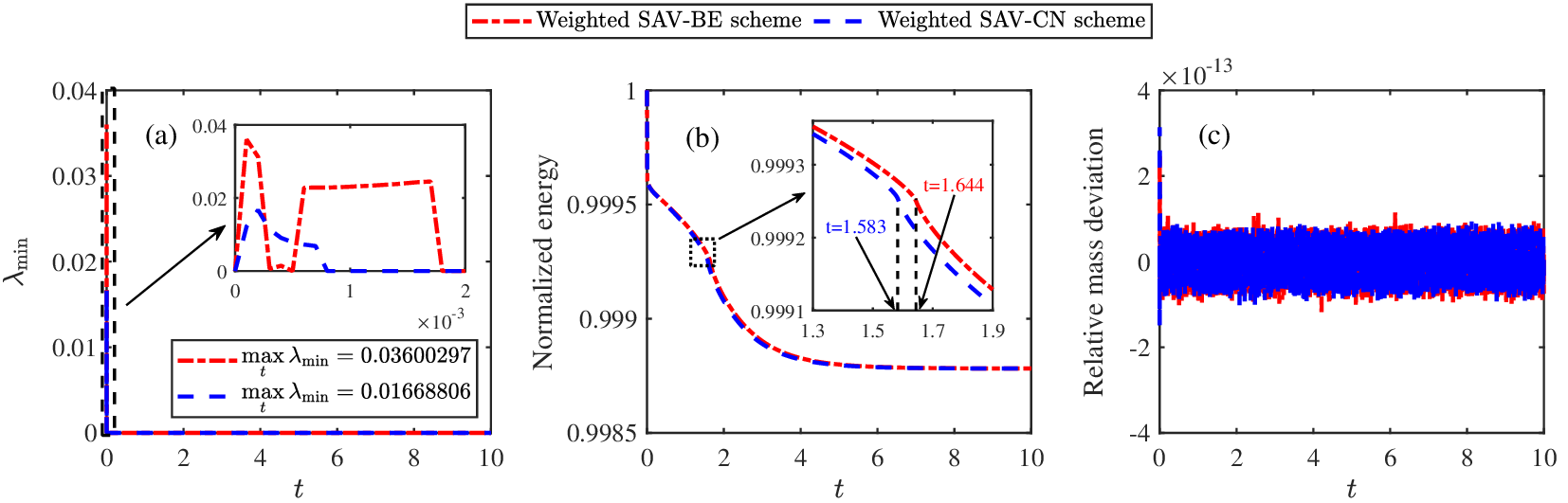}
\caption{The evolution of (a) weighting coefficient $\lambda_{\min}$, (b) normalized energy, and (c) relative mass deviation
         simulated using the weighted SAV-BE scheme \eqref{eqc3} and weighted SAV-CN scheme \eqref{eqd1}, respectively.}
\label{figt4}
\end{figure}

Fig.~\ref{figc4} depicts the evolution process under the second-order weighted SAV-CN scheme.
It can be observed that the central hole in the torus gradually shrinks and eventually disappears,
leading the system to evolve into a sphere and reach equilibrium state.
At the phase transition (i.e., at $t=1.583$), the energy curve still exhibits significant changes.
Next, in Fig.~\ref{figt4}, we compare the first-order and second-order weighted SAV schemes.
It is evident that the second-order scheme indicates a smaller weighting coefficients,
while the evolution process of the first-order scheme is slightly slower.
Furthermore, the both schemes successfully maintain the conservation properties of mass.

\section{Conclusion}
\label{sec5}

By integrating a weighted concept,
we consolidate the methodologies of the nonlinear energy-based SAV and Lagrange multiplier-based SAV within a unified framework.
Our approach notably improves the discrete energy approximation,
making it closer to the original energy compared to the nonlinear energy-based approach,
and effectively avoids potential unsolvability issues associated with the Lagrange multiplier approach.
The numerical results robustly validate that the weighted SAV-BE and SAV-CN schemes exhibit first-order and second-order accuracy in time, respectively, while preserving energy stability. 
Due to the wide application of SAV-like approaches in addressing various physical dissipative problems in scientific computing,
our future works intend to extending the weighted concept to tackle more challenges, such as discrepancies
between discrete and original energies under large time steps, as well as other unresolved issues associated with the Lagrange multiplier approach.

\begin{acknowledgement}
The numerical calculations in this paper have been done on the supercomputing system in the Supercomputing Center of Wuhan University.
\end{acknowledgement}
\medskip
\noindent{\textbf{Funding}} ~ This work was partially supported by the National Natural Science Foundation of China
(Nos. 12001210, 12131010, 12271414, 12125103, 12071362, 12301558).

\medskip
\noindent{\textbf{Data Availability}} ~The code for the current study is available from the corresponding author on reasonable
request.

\section*{Declarations}
{\textbf{Conflict of interest}} ~The authors declare no competing interests.



\end{document}